\input amstex
\input epsf
\documentstyle{amsppt}
\nologo\footline={}

\hsize450pt

\def\PU{\mathop{\text{\rm PU}}}
\def\Rep{\mathop{\text{\smc Rep}}}
\def\tr{\mathop{\text{\rm tr}}}
\def\SU{\mathop{\text{\rm SU}}}
\def\B{\mathop{\text{\rm B}}}
\def\dist{\mathop{\text{\rm dist}}}
\def\ta{\mathop{\text{\rm ta}}}
\def\Arg{\mathop{\text{\rm Arg}}}
\def\area{\mathop{\text{\rm area}}}
\def\Im{\mathop{\text{\rm Im}}}
\def\U{\mathop{\text{\rm U}}}
\def\rk{\mathop{\text{\rm rk}}}
\def\Re{\mathop{\text{\rm Re}}}
\def\SL{\mathop{\text{\rm SL}}}

\topmatter
\title
Hyperbolic $2$-spheres with cone singularities
\endtitle
\author
Sasha Anan$'$in, Carlos H.~Grossi, Jaejeong Lee, Jo\~ao dos Reis jr.
\endauthor
\address
Departamento de Matem\'atica, ICMC, Universidade de S\~ao Paulo, Caixa Postal 668,\newline
13560-970--S\~ao Carlos--SP, Brasil
\endaddress
\email
sasha\_a\@icmc.usp.br
\endemail
\address
Departamento de Matem\'atica, ICMC, Universidade de S\~ao Paulo, Caixa Postal 668,\newline
13560-970--S\~ao Carlos--SP, Brasil
\endaddress
\email
grossi\@icmc.usp.br
\endemail
\address
Korea Institute for Advanced Study (KIAS)\newline85 Hoegiro Dongdaemun-gu, Seoul 02455, Republic of Korea
\endaddress
\email
jjlee\@kias.re.kr
\endemail
\address
Departamento de Matem\'atica, ICMC, Universidade de S\~ao Paulo, Caixa Postal 668,\newline
13560-970--S\~ao Carlos--SP, Brasil
\endaddress
\email
joao.reis.reis\@usp.br
\endemail
\abstract
We study the space $C(a_0,a_1,\dots,a_n)$ of hyperbolic $2$-spheres with cone points of prescribed apex curvatures
$2a_0,2a_1,\dots,2a_n\in]0,2\pi[$ and some related spaces. For $n=3$, we get a detailed description of such spaces. The euclidean
$2$-spheres were considered by W.~P.~Thurston: for $n=4$, the corresponding spaces provide the famous $7$ examples of
nonarithmetic compact holomorphic $2$-ball quotients previously constructed by Deligne-Mostow.
\endabstract
\endtopmatter
\document

{\hfill\it To Misha Kapovich, a geometric geometer}

\bigskip

\centerline{\bf1.~Introduction}

\medskip

In the seminal work [Thu], W.~P.~Thurston studied the space of flat $2$-spheres with cone points of prescribed apex curvatures.
Such a space possesses a structure of complex hyperbolic manifold, though, is not in general complete. Its completion has a
natural structure of complex hyperbolic cone manifold which turns out to be an orbifold exactly when satisfies a certain simple
{\sl orbifold condition\/} on the prescribed apex curvatures. In this way, Thurston reobtained the celebrated 7 nonarithmetic
compact holomorphic $2$-ball quotients (for short, nonarithmetic {\it quotients\/}) that were constructed by Deligne-Mostow
[DM1], [DM2] about 30--35 years ago. As we are not going to distinguish a quotient $Q$ from a unramified finite cover of $Q$, we
assume $Q$ to be a manifold. Recently, 4 more nonarithmetic quotients were constructed in~[DPP]. The latter 4 are quite similar
to the former 7. Each of 11 contains at least four {\it$\Bbb C$-fuchsian curves,} i.e., smooth holomorphic totally geodesic
compact curves.

In turn, the arithmetic quotients are known to be of 2 types [McR]: the {\it first\/} type and the {\it second\/} type which is related to division algebras. Taking into account the Mostow-Prasad rigidity of quotients, there are countably many
arithmetic quotients of either type. The types can be distinguished by the presence of a $\Bbb C$-fuchsian curve: a quotient of
the first type possesses such a curve, whereas that of the second one does not. Using $\Bbb C$-fuchsian curves, we define the
type of an arbitrary quotient with the same criterion and arrive at the following questions.

\medskip

{\bf1.1.~Problem.} Is there any nonarithmetic quotient of the second type? Do there exist infinitely many nonarithmetic quotients
of the first/second type?

\medskip

{\bf1.2.~Outline.} In this paper, we develop a few tools allowing to study the space $C(a_0,a_1,\dots,a_n)$ of $2$-spheres with
cone points of prescribed apex curvatures $2a_0,2a_1,\dots,2a_n$, where the set of the other points is endowed with the {\sl
spherical\/} geometry of curvature $\sigma:=1$ or with the {\sl hyperbolic\/} geometry of curvature $\sigma:=-1$. Hopefully,
$C(a_0,a_1,\dots,a_n)$ or related spaces can shed some light on the above problem.

\medskip

{\bf1.2.1.~Toy example.} In order to clarify the general picture, we first glance at the toy case of $4$ cone points of the fixed
apex curvatures $2a,\pi,\pi,\pi$, where $a\in]0,\frac\pi2[$ when $\sigma=1$ and $a\in]\frac\pi2,\pi[$ when $\sigma=-1$. In this
particular situation, each spherical or hyperbolic $2$-sphere can be cut into a solid triangle (closed disc) bounded by a
triangle and living in the model space, i.e., in the round sphere or in the hyperbolic plane. Indeed, we simply join the cone
point $p$ of apex curvature $2a$ with the other 3 cone points by means of simple geodesic segments so that the segments intersect
only in $p$ (say, taking the shortest ones) and then cut the $2$-sphere along the segments. The point $p$ turns into the vertices
of the triangle and the other 3 cone points, into the middle points of the sides of the triangle. By splitting the solid triangle
along a median and then gluing the corresponding half-sides, we {\it modify\/} the original triangle with respect to the median.
As it happens, the $2$-spheres corresponding to triangles $t_1$ and $t_2$ are isometric (orientation preserved) iff $t_1$ and
$t_2$ differ by finitely many such modifications. This allows us to describe the space $C(a,\frac\pi2,\frac\pi2,\frac\pi2)$ as
the quotient of a topological disc $S$, the space of counterclockwise oriented triangles of the fixed area $\sigma(\pi-2a)$, by
the group $M_3$ generated by $3$ involutions which are nothing but the modifications acting on $S$. In a certain sense, this
construction is reminiscent of [WaW].

\medskip

{\bf1.2.2.~General case.} Denote by $G$ the Lie group of orientation-preserving (= holomorphic) isometries of the model space.
So, $G:=\PU(2)$ when $\sigma=1$ and $G:=\PU(1,1)$ when $\sigma=-1$. Let $c_0,c_1,\dots,c_n$ stand for the cone points of some
$2$-sphere $\Sigma\in C(a_0,a_1,\dots,a_n)$. Then we get a holonomy representation $\varrho:F_n\to G$, where
$F_n:=\pi_1\big(\Sigma\setminus\{c_0,c_1,\dots,c_n\}\big)$ is the free group of rank $n$. In this way, we~obtain a well-defined
holonomy map $h:C(a_0,a_1,\dots,a_n)\to\Rep(F_n,G)/G$, $\Sigma\mapsto[\varrho]$, to the space of representations $\Rep(F_n,G)$
considered modulo conjugation by $G$.

As in 1.2.1, we can join $c_0$ with $c_1,\dots,c_n$ by means of $n$ simple geodesic segments so that the segments intersect only
in $c_0$ (say, taking the shortest ones) and then cut $\Sigma$ along the segments thus getting a $2n$-gon $P$ that bounds a
topological closed disc not necessarily embeddable into the model space. Labeling the odd vertices of $P$ that originate from
$c_1,\dots,c_n$ with the same symbols and the even ones, that originate from $c_0$, with $p_1,\dots,p_n$, we see that the interior
angles of $P$ at the $p_j$'s sum to $2\pi-2a_0$, that the interior angle of $P$ at $c_j$ equals $2\pi-2a_j$, and that the sides
of $P$ adjacent at $c_j$ have equal length for all $j=1,\dots,n$. In some sense, $P$ bounds a fundamental domain for $\Sigma$.

A choice of the above segments distinguishes some generators $r_0,r_1,\dots,r_n$ of $F_n$ such that
$F_n=\langle r_0,r_1,\dots,r_n\mid r_n\dots r_1r_0=1\rangle$. A suitable braid group $\overline M_n$ (see Definition 4.2) acts on
$F_n$ and preserves the set of conjugacy classes of the $r_j$'s. Another choice of the segments is given at the level of
representations by an element $m\in\overline M_n$ that transforms $[\varrho]$ into $[\varrho\circ m]$. (In general, such elements
$m$ do not form a subgroup.) Hence, the image of the holonomy map $h$ lives in a sort of relative character variety
$R(a_0,a_1,\dots,a_n)\subset\Rep(F_n,G)/G$ formed by all the representations $[\varrho]$ such that the conjugacy classes of
$\varrho r_0,\varrho r_1,\dots,\varrho r_n$ are the counterclockwise rotations in the model space by the angles
$2a_0,2a_1,\dots,2a_n$, listed perhaps in different order. In alternative words, $R(a_0,a_1,\dots,a_n)$ is the space of relations
between rotations in prescribed conjugacy classes (the order of the classes is not fixed).

\medskip

{\bf1.2.3.~Convex hyperbolic $2n$-gons.} Let us look more closely at the case of $\sigma=-1$ with
$a_0,a_1,\dots,a_n\in[\frac\pi2,\pi[$ such that $\sum_ja_j>2\pi$. Here, any $2n$-gon $P$ is convex, i.e., its interior angles
are all in $]0,\pi]$, hence, the corresponding closed disc is embeddable into the model space. In this case, the space
$C(a_0,a_1,\dots,a_n)$ of hyperbolic 2-spheres with cone points of apex curvatures $2a_0,2a_1,\dots,2a_n$ can be described as
the quotient of a component of $R(a_0,a_1,\dots,a_n)$ by the action of $\overline M_n$ (see Section~5 and Theorem 5.1; see also
Section 1.3).

\medskip

{\bf1.2.4.~Geometry on $C(a_0,a_1,\dots,a_n)$ and on the related spaces.} Since we hope to endow $C(a_0,a_1,\dots,a_n)$ with a
structure of complex hyperbolic orbifold, it is natural to conjecture that such a geometry descents from a component $C$ of
$R(a_0,a_1,\dots,a_n)$ and that $\overline M_n$ acts discretely on $C$ by isometries. (There is a more subtle approach to
constructing spaces with geometry in a similar vein, but let us stick first to a simplest variant.)

Actually, even in the flat case considered by Thurston, there may exist no preferred hyperbolic geometry on $C$. (Well, in the
flat case, one can always choose the structure related to the oriented area function; it does provide a complex hyperbolic
structure.) The choice of a `good' hyperbolic structure on $C$ is the principal challenge in our project, and the first step is
to describe the components of $R(a_0,a_1,\dots,a_n)$ where $\overline M_n$ acts discretely.

One can establish the lack of geometry of constant curvature on a component of $R(a_0,a_1,a_2,a_3)$, compatible with the action
of $\overline M_3$, if this action is not discrete.

\medskip

{\bf1.2.5.~Brief sketch of the exposition.} Every representation $[\varrho]\in R(a_0,a_1,\dots,a_n)$ can be inter\-preted as a
labeled $2n$-gon, a closed piecewise geodesic oriented path, whose consecutive vertices $c_1,p_1,\dots,c_n,p_n:=c_0$ are given by
$p_j:=R_jp_{j-1}$ for all $1\leqslant j\leqslant n$ (the indices are modulo $n$), where $c_j$ stands for the fixed point of the
rotation $R_j:=\varrho r_j$, $0\leqslant j\leqslant n$, and, without loss of generality, $R_0$~is a counterclockwise rotation by
$2a_0$.

It is convenient to deal with the $2n$-gons with a forgotten label $c_0$. In other words, we consider the quotient
$S(a_0,a_1,\dots,a_n)$ of $R(a_0,a_1,\dots,a_n)$ by a suitable order $n$ cyclic subgroup in $\overline M_n$. Thus, we are allowed
to place the fixed point of the rotation $R_0$ at any vertex $p_j$. In terms of relations, the generator of the cyclic group
transforms the relation $R_n\dots R_1R_0=1$ into the relation $R_1R_n\dots(R_1R_0R_1^{-1})=1$. A certain braid group $M_n$ (see
Definition 4.3) generated by simple modifications (surgeries) of the $2n$-gons acts on $S(a_0,a_1,\dots,a_n)$ so that the map
$R(a_0,a_1,\dots,a_n)\to S(a_0,a_1,\dots,a_n)$ induces a bijection between the $\overline M_n$-orbits and the $M_n$-orbits.
Clearly, the discreteness of $\overline M_n$ is equivalent to that of $M_n$.

Next step is to describe the components of $S(a_0,a_1,\dots,a_n)$. The cyclic order of $a_1,\dots,a_n$ in $P$ defines the {\it
type\/} of a $2n$-gon $P$ and $M_n$ acts transitively on the types. Thus, we need to describe the components of a given type
and to find those components where $M$, the index $(n-1)!$ subgroup of the type-preserving elements of $M_n$, acts discretely.

We proceed by induction on $n$. Without loss of generality, we assume $R_0$ and $R_n$ to be counterclockwise rotations by $2a_0$
and $2a_n$. The conjugacy class of $R_0R_n$ is given exactly by the distance $d_n$ between $c_n$ and $p_n$. Therefore, by
induction, we know the topology of a fibre of the map $d_n:S(a_0,a_1,\dots,a_n)\to[0,\infty[$. Indeed, if $R_0R_n$ is a
counterclockwise rotation by $2a$, $a\in]0,\pi[$, then the fibre in question is $S(a,a_1,\dots,a_{n-1})\times C_{R_0R_n}$, where
$C_{R_0R_n}$ stands for the centralizer of $R_0R_n$ in~$G$, topologically a circle. If $R_0R_n=1$ (hence, $d_n=0$), then the
fibre equals $S(a_1,\dots,a_{n-1})$. If~$R_0R_n$ is a parabolic or hyperbolic isometry, then the fibre is
$S(a^*,a_1,\dots,a_{n-1})\times C_{R_0R_n}$, where $a^*$ stands for the conjugacy class of $R_0R_n$ and $C_{R_0R_n}$ is
topologically a line. (This means that, for the sake of induction, we need to consider the relations between isometries in
prescribed conjugacy classes, not~necessarily ellip\-tic ones. We do not do so in this particular paper because we study here
mostly the case of $n=3$ and $\sigma=-1$.) Gathering information about the topology of the fibres of $d_n$, including that at the
singular points of $d_n$, we arrive at the description of the components and their topology (see~Proposition~4.5.2).

Also, in the case $n=3$ and $\sigma=-1$, we visualize 3 curves $A_1,A_2,A_3$ that are respectively fixed point sets of the
involutions $n_1,n_2,n_3$ that generate the group $M_3$ (here we identify the components of both types). We show that $M_3$ does
not act discretely on a component if the component is not smooth (see Remark 4.6.4) or if a couple of the above curves intersects
in the component. Finally, we~describe all components where $M_3$ acts discretely (see Theorem 4.6.3) : for any
$a_0,a_1,a_2,a_3\in]0,\pi[$, unless $a_0=a_1=a_2=a_3=\frac\pi2$ (the case of nonsmooth $S(a_0,a_1,a_2,a_3)$), there exists
a component (of~a~given type) where $M_3$ acts discretely iff $a_j+a_k\leqslant\pi$ for all $j\ne k$ or $\pi\leqslant a_j+a_k$
for all $j\ne k$, and~such a component is unique. Taking an ideal triangle on the hyperbolic plane and the group generated by the
reflections in its sides $A_1,A_2,A_3$, we get an adequate topological picture of the action of $M_3$ on a component in question.

It is worthwhile mentioning that any component of $S(a_0,a_1,a_2,a_3)$ is a surface in $\Bbb R^3$ given by one of the equations
(3.5) or (3.7).

\medskip

{\bf1.3.~Related works.} Being back to the toy example 1.2.1, hence, taking in the equations (3.5) and (3.7)
$a_1=a_2=a_3:=\frac\pi2$, and changing the variables with respect to $x_j:=\pm2s(2t_j-1)$, we arrive at the equation
$x_1^2+x_2^2+x_3^2-x_1x_2x_3-2=t$, where $s:=\sin a_0$, $t:=-\tr R_{c,k}$, $k:=e^{a_0i}$, and $R_{c,k}\in\SU(1,1)$ (see the
definition in the beginning of Section 2) represents the counterclockwise rotation by $2a_0$ about $c$ in the hyperbolic plane.
In other words, our group $M_3$ in this case is commensurable with the group $\Gamma$ from [Gol] in the case $|t|<2$. If we would
admit in the toy example a parabolic or hyperbolic isometry in place of $R_{c,k}$, we could deal in fact with the general case
studied by W.~M.~Goldman [Gol] (the~parabolic isometry would correspond to the case studied by P.~Waterman and S.~Wolpert [WaW]).

The condition $\sum_{j=0}^na_j>2\pi$ (or $\sum_{j=0}^na_j<2\pi$) on the prescribed apex curvatures, necessary for the existence
of a hyperbolic (respectively, spherical) $2$-sphere $\Sigma$ with the indicated cone singularities, is clearly sufficient.
Moreover,
A.~D.~Alexandrov proved that, for any such $\Sigma$, there exists a unique convex compact polyhedron $P$ in the real hyperbolic
$3$-space of curvature $-1$ (respectively, in the round $3$-sphere of curvature $1$) whose boundary $\partial P$ is isometric to
$\Sigma$ with respect to the inner metric on $\partial P$ (see [Ale,~3.6.4, p.~190] for the uniqueness and [Ale, 5.3.1, p.~261]
for the existence).

G.~Mondello and D.~Panov [MPa] found necessary and almost sufficient conditions for the existence of a spherical $2$-sphere with
cone points of given cone angles, each of the angles is allowed to exceed $2\pi$. The corresponding space can be quite useful
when dealing with Problem 1.1. In the context of convex polyhedra in real hyperbolic $3$-space, C.~D.~Hodgson and I.~Rivin [HRi]
described all spherical $2$-spheres with all cone angles $>2\pi$ and all closed geodesics of length $>2\pi$.

\medskip

{\bf Acknowledgements.} We are very grateful to KIAS and to professor Sungwoon Kim for organizing `Workshop on geometric
structures, Hitchin components, and representation varieties' and supporting our collaboration (see
http://home.kias.re.kr/MKG/h/WGS2015/). We also thank the referee for remarks that improved our exposition. The first three authors were supported by FAPESP, grant 2015/25809-2. The forth author was partially supported by FAPESP, grant
2014/26295-0.

\bigskip

\centerline{\bf2.~Preliminary remarks}

\medskip

We are going to use the following settings and notation similar to those in [AGr].

Let $\sigma=\pm1$. In what follows, $V$ denotes a $2$-dimensional $\Bbb C$-linear space equipped with a hermitian form
$\langle-,-\rangle$ of signature $\big(1+\frac{\sigma+1}2,\frac{\sigma-1}2\big)$. Denote
$\B V:=\{p\in\Bbb P_\Bbb CV\mid\langle p,p\rangle>0\}$. Then $\B V$ is the model space of curvature $\sigma$, a round sphere if
$\sigma=1$ and a hyperbolic disc if $\sigma=-1$. In~fact, when $\sigma=-1$, the {\it Riemann-Poincar\'e sphere\/}
$\Bbb P_\Bbb CV$ is glued from two hyperbolic discs $\B V$ and $\B'V:=\{p\in\Bbb P_\Bbb CV\mid\langle p,p\rangle<\nomathbreak0\}$
along the {\it absolute\/} $\partial\B V:=\{p\in\Bbb P_\Bbb CV\mid\langle p,p\rangle=0\}$.

The distance $\dist(p_1,p_2)$ between points $p_1,p_2\in\B V$ is given by the formulae
$\cos^2\frac{\dist(p_1,p_2)}2=\ta(p_1,p_2)$ for $\sigma=1$ and $\cosh^2\frac{\dist(p_1,p_2)}2=\ta(p_1,p_2)$ for $\sigma=-1$. So,
the distance is a monotonic function of the {\it tance\/} defined as
$\ta(p_1,p_2):=\frac{\langle p_1,p_2\rangle\langle p_2,p_1\rangle}{\langle p_1,p_1\rangle\langle p_2,p_2\rangle}$.

For $z\in\Bbb C\setminus0$, we denote by $\arg z$ and $\Arg z$ the argument functions taking values in $]{-}\pi,\pi]$ and in
$[0,2\pi[$, respectively.

\medskip

{\bf2.1.~Remark {\rm[AGr, Example 6.1]}.} {\sl The oriented area of a triangle with pairwise nonorthogonal vertices\/
$p_1,p_2,p_3\in\B V$ equals\/ $\area(p_1,p_2,p_3)=-2\sigma\arg(g_{12}g_{23}g_{31})$, where\/ $[g_{jk}]$ stands for the Gram
matrix of the\/ $p_j$'s. {\rm(}If\/~$\sigma=-1$, it is possible to observe that\/ $\arg$ takes values in\/
$]{-}\frac\pi2,\frac\pi2[$ when calculating\/ area by the above formula\/{\rm;} if\/ $\sigma=1$, then\/ $\arg$ can take any
value.{\rm)} In particular, the triangle is counterclockwise oriented iff\/ $-\sigma\Im(g_{12}g_{23}g_{31})>0$.}

\medskip

For $c\in\B V$ and $k\in\Bbb C$ such that $|k|=1$, the rule
$R_{c,k}:p\mapsto(\overline k-k)\frac{\langle p,c\rangle}{\langle c,c\rangle} c+kp$ defines a $\Bbb C$-linear map
$R_{c,k}:V\to V$.

\medskip

{\bf2.2.~Remark.} {\sl Let\/ $c,p\in\B V$ and let\/ $k,k'\in\Bbb C$ be such that\/ $|k|=|k'|=1$. Then\/ $R_{c,k}$
provides a rotation in\/ $\B V$. More precisely, $R_{c,k}$ is the rotation by\/ $\Arg k^2$ in the counterclockwise sense
about\/~$c$. Moreover, $R_{c,k}\in\SU V$, $R_{c,-k}=-R_{c,k}$, $R_{c,k}R_{c,k'}=R_{c,kk'}$,
$\langle c,R_{c,k}p\rangle=k\langle c,p\rangle$, and\/
$\langle R_{c,k}p,p\rangle=\big(\ta(p,c)(\overline k-k)+k\big)\langle p,p\rangle$. If\/ $\sigma=1$ and $\langle c,c'\rangle=0$,
then\/ $R_{c,k}R_{c',k}=1$.}

\medskip

{\bf Proof.} Let $c'\in\Bbb P_\Bbb CV$ be the point orthogonal to $c$ and take representatives such that $\langle c,c\rangle=1$
and $\langle c',c'\rangle=\sigma$. Every point $x\in\B V$ can be written in the form $x=c+\eta c'$ for some unique
$\eta\in\Bbb C$. When $\sigma=-1$, we have $|\eta|<1$ and this corresponds to identifying $\B V$ with the unitary disc in
$\Bbb C$ centred at the origin; when $\sigma=1$, we identified $\Bbb P_\Bbb CV\setminus c'$ with $\Bbb C$ (in both cases, the
point $c$ is mapped to the origin). Under this identification, $R_{c,k}$ sends the point $x$ to $c+k^2\eta c'$, i.e., it is
nothing but multiplication by the unitary complex number $k^2$. In the orthonormal
basis $c,c'$, the map $R_{c,k}$ has the form $\left[\smallmatrix\overline k&0\\0&k\endsmallmatrix\right]$, so it belongs to
$\SU V$. The equalities $R_{c,-k}=-R_{c,k}$, $R_{c,k}R_{c,k'}=R_{c,kk'}$, $\langle c,R_{c,k}p\rangle=k\langle c,p\rangle$, and
$\langle R_{c,k}p,p\rangle=\big(\ta(p,c)(\overline k-k)+k\big)\langle p,p\rangle$ follow from straightforward calculations.
Finally, in order to show that $R_{c,k}R_{c',k}=1$, it suffices to note that $R_{c,k}R_{c',k}c=c$ and $R_{c,k}R_{c',k}c'=c'$.


\medskip

Let $\gamma$ be a geodesic in $\B V$. We denote by $I_\gamma$ the (antiholomorphic) reflection in $\gamma$. Given a couple of
geodesics $\gamma,\gamma'$ intersecting at $c\in\B V$, denote by $\alpha_c(\gamma,\gamma')\in[0,\pi[$ the {\it oriented angle\/}
at $c$ from $\gamma$ to $\gamma'$ counted in the counterclockwise sense.

\medskip

{\bf2.3.~Remark.} {\sl Let\/ $1\ne R\in\PU V$ be a nontrivial orientation-preserving isometry of\/ $\B V$. Then there exist
distinct geodesics $\gamma,\gamma'$ such that\/ $R=I_{\gamma'}I_\gamma$ in\/ $\PU V$.

The geodesics are ultraparallel, tangent, or intersecting iff\/ $R$ is respectively hyperbolic, parabolic, or~elliptic. The
couple of geodesics is unique up to the action of the centralizer\/ $C_R$ of\/ $R$ in\/ $\PU V$ and the action of the centralizer
is free in the hyperbolic and parabolic cases. In the elliptic case, the stabilizer of the couple is given by the reflection in a
point of the intersection of the geodesics, i.e., by\/ $R_{c,i}$, where\/ $c$ is a point of the intersection. In this case,
$I_{\gamma'}I_\gamma=R_{c,k}\in\PU V$, where\/ $k:=e^{i\alpha_c(\gamma,\gamma')}$.

Topologically, $C_R$ is a circle if\/ $R$ is elliptic and a line, if\/ $R$ is hyperbolic or parabolic.}

\medskip

{\bf2.4.~Remark.} {\sl Let\/ $\sigma=-1$. For any\/ $a_1,a_2\in]0,\pi[$ such that\/ $a_1+a_2\ne\pi$, there exists unique\/
$t_0(a_1,a_2)>1$ satisfying the following property. Given distinct points\/ $q_1,q_2\in\B V$, denote by\/ $\gamma_1,\gamma_2$ the
geodesics such that\/ $q_j\in\gamma_j$, $\alpha(\gamma,\gamma_1)=a_1$, and\/ $\alpha(\gamma_2,\gamma)=a_2$, where\/ $\gamma$
stands for the geodesic joining\/ $q_1,q_2$. The geodesics\/ $\gamma_1,\gamma_2$ intersect, are tangent, or are ultraparallel
if\/ $\ta(q_1,q_2)<t_0(a_1,a_2)$, $\ta(q_1,q_2)=t_0(a_1,a_2)$, or $\ta(q_1,q_2)>t_0(a_1,a_2)$, respectively. In the case\/
$a_1+a_2=\pi$, the geodesics\/ $\gamma_1,\gamma_2$ are always ultraparallel and we put\/ $t_0(a_1,a_2):=1$ in this case.}

\medskip

{\bf2.5.~Remark.} {\sl Let\/ $\gamma_1,\gamma_2$ be ultraparallel geodesics and let\/ $a_1,a_2\in]0,\pi[$. Then there exists a
unique geodesic\/ $\gamma$, depending smoothly on\/ $\gamma_1,\gamma_2$, such that\/ $\alpha(\gamma_1,\gamma)=a_1$ and\/
$\alpha(\gamma,\gamma_2)=a_2$.}

\medskip

{\bf Proof.} Let $v_1,v'_1$ and $v_2,v'_2$ be vertices of $\gamma_1$ and of $\gamma_2$ such that the triangles $(v'_1,v_1,v_2)$
and $(v'_1,v_2,v'_2)$ are both clockwise oriented. Then the angle from $[v'_1,c]$ to $[c,v_2]$ grows monotonically from $0$ to
$\pi$ while $c\in\gamma_1$ runs from $v_1$ to $v'_1$. At some $c\in\gamma_1$, this angle equals $a_1$. Similarly, there exists
$c'\in[c,v'_1]$ such that the angle from $[v'_1,c']$ to $[c',v'_2]$ equals $a_1$. It follows that, for any $p\in[c,c']$, the~
geodesic $\gamma_p$ that passes through $p$ and such that $\alpha(\gamma_1,\gamma_p)=a_1$ intersects $\gamma_2$. For $p=c$, this
intersection is $v_2$ and, for $p=c'$, it is $v'_2$. Hence, the angle $\alpha(\gamma_p,\gamma_2)$ varies from $0$ to $\pi$. By
continuity, we obtain the existence of the desired $\gamma$. If we have another $\gamma'$, then $\gamma$ and $\gamma'$ are
ultraparallel. So,~we~get a simple quadrangle whose interior angles sum to $2\pi$, a contradiction.

\medskip

{\bf2.6.~Remark.} {\sl Let\/ $\sigma=-1$ and\/ $a_0,a_1,a_2,a_3\in]0,\pi[$. Then, for any\/ $t>t_0(a_0,a_3)$, there exist
geometrically unique geodesics\/ $\gamma,\gamma_0,\gamma_1,\gamma_2$ such that\/ $\alpha(\gamma_0,\gamma_1)=a_0$,
$\alpha(\gamma_1,\gamma)=a_1$, $\alpha(\gamma,\gamma_2)=a_2$, $\alpha(\gamma_2,\gamma_0)=a_3$, and\/ $\ta(c_3,c_0)=t$
{\rm(}hence, $\gamma_1,\gamma_2$ are ultraparallel\/{\rm)}, where\/ $c_0$ and\/ $c_3$ stand for the intersection points of\/
$\gamma_0,\gamma_1$ and of\/ $\gamma_0,\gamma_2$, respectively. Moreover, the four geodesics depend smoothly on\/
$t>t_0(a_0,a_3)$.}

\medskip

{\bf Proof.} By Remark 2.4, taking points $c_0,c_3$ on the tance $\ta(c_0,c_3)=t$, we get ultraparallel geodesics
$\gamma_1,\gamma_2$ such that $c_0\in\gamma_1$ with $\alpha(\gamma_0,\gamma_1)=a_0$ and $c_3\in\gamma_2$ with
$\alpha(\gamma_2,\gamma_0)=a_3$, where $\gamma_0$ stands for the geodesic joining $c_0,c_3$. By Remark 2.5, there exists a unique
geodesic $\gamma$ such that $\alpha(\gamma_1,\gamma)=a_1$ and $\alpha(\gamma,\gamma_2)=a_2$.

\medskip

{\bf2.7.~Remark.} {\sl Let\/ $\sigma=-1$. Given distinct geodesics\/ $\gamma_1,\gamma_2$ intersecting at some point\/ $c\in\B V$,
denote\/ $a:=\alpha(\gamma_2,\gamma_1)$ and pick some\/ $a_1,a_2\in]0,\pi[$. Then there exists a couple of geodesics\/
$\gamma,\gamma'$ {\rm(}subject to $\gamma'=R_{c,i}\gamma${\rm)} such that\/
$\alpha(\gamma_1,\gamma)=\alpha(\gamma_1,\gamma')=a_1$ and\/ $\alpha(\gamma,\gamma_2)=\alpha(\gamma',\gamma_2)=a_2$ iff

\smallskip

$\bullet$ $a+a_1+a_2\leqslant\pi$, when the triangle\/ $(c,c_1,c_2)$ is clockwise oriented or degenerate\/ {\rm(}i.e.,
$c=c_1=c_2${\rm)}, or

$\bullet$ $2\pi\leqslant a+a_1+a_2$, when the triangle\/ $(c,c_1,c_2)$ is counterclockwise oriented or degenerate,

\smallskip

\noindent
where\/ $c_j$ stands for the intersection point of\/ $\gamma$ and\/ $\gamma_j$. Such geodesics\/ $\gamma,\gamma'$ are unique when
exist and depend continuously on\/ $\gamma_1,\gamma_2$, and\/ $a$.

The condition\/ $\gamma=\gamma'$ is equivalent to\/ $a+a_1+a_2=\pi$ or\/ $a+a_1+a_2=2\pi${\rm;} equivalently, this means that\/
$c=c_1=c_2$.}

\medskip

{\bf2.8.~Remark.} {\sl Let\/ $\sigma=-1$. Suppose that\/ $a_0,a_1,a_2,a_3\in]0,\pi[$ satisfy\/ $a_0+a_1+a_2+a_3<\pi$ or\/
$3\pi<a_0+a_1+a_2+a_3$. Then there exist geodesics\/ $\gamma,\gamma_0,\gamma_1,\gamma_2$ such that\/ $\gamma_1,\gamma_2$
intersect, $\alpha(\gamma_0,\gamma_1)=a_0$, $\alpha(\gamma_1,\gamma)=a_1$, $\alpha(\gamma,\gamma_2)=a_2$, and\/
$\alpha(\gamma_2,\gamma_0)=a_3$.}

\medskip

{\bf Proof.} If $a_0+a_1+a_2+a_3<\pi$, then $0<a_0+a_3<\pi-a_1-a_2<\pi$. So, we can pick $a\in]0,\pi[$ such that
$a_0+a_3<a<\pi-a_1-a_2$. Hence, there is a counterclockwise oriented triangle $(c,c_0,c_3)$ with the interior angles
$\pi-a,a_0,a_3$, respectively. Denote by $\gamma_0,\gamma_1,\gamma_2$ the geodesics such that $c_0,c_3\in\gamma_0$,
$c,c_0\in\gamma_1$, and $c_3,c\in\gamma_2$. As $a=\alpha(\gamma_2,\gamma_1)$ and $a+a_1+a_2<\pi$, by Remark 2.7, there exists a
geodesic $\gamma$ with $\alpha(\gamma_1,\gamma)=a_1$ and $\alpha(\gamma,\gamma_2)=a_2$.

If $3\pi<a_0+a_1+a_2+a_3$, then $0<2\pi-a_1-a_2<a_0+a_3-\pi<\pi$. So, we can pick $a\in]0,\pi[$ such that
$2\pi-a_1-a_2<a<a_0+a_3-\pi$. Hence, there is a clockwise oriented triangle $(c,c_0,c_3)$ with the interior angles
$a,\pi-a_0,\pi-a_3$, respectively. Denote by $\gamma_0,\gamma_1,\gamma_2$ the geodesics such that $c_0,c_3\in\gamma_0$,
$c,c_0\in\gamma_1$, and $c_3,c\in\gamma_2$. As $a=\alpha(\gamma_2,\gamma_1)$ and $2\pi<a+a_1+a_2$, by Remark 2.7, there exists a
geodesic $\gamma$ with $\alpha(\gamma_1,\gamma)=a_1$ and $\alpha(\gamma,\gamma_2)=a_2$.

\medskip

{\bf2.9.~Remark.} {\sl Let\/ $\gamma_1,\gamma_2$ be distinct tangent geodesics and let\/ $a_1,a_2\in]0,\pi[$. Denote by\/ $v_1,v$
and by\/ $v,v_2$ the vertices of\/ $\gamma_1$ and of\/ $\gamma_2$, respectively. Then there exists a geodesic\/ $\gamma$ such
that\/ $\alpha(\gamma_1,\gamma)=a_1$ and\/ $\alpha(\gamma,\gamma_2)=a_2$ iff\/ $a_1+a_2<\pi$ and the triangle\/ $(v_1,v,v_2)$ is
counterclockwise oriented or\/ $\pi<a_1+a_2$ and the triangle\/ $(v_1,v,v_2)$ is clockwise oriented. Such a geodesic\/ $\gamma$
is unique if exists.

In particular, any sufficiently small deformation of\/ $\gamma_1,\gamma_2$ still allows a geodesic\/ {\rm(}a couple of
geodesics{\rm)} with the required angles.}

\medskip

{\bf2.10.~Remark.} {\sl Let\/ $\sigma=-1$. Suppose that\/ $a_0,a_1,a_2,a_3\in]0,\pi[$ satisfy\/ $a_0+a_3<\pi<a_1+a_2$ or\/
$a_1+a_2<\pi<a_0+a_3$. Then there exist geodesics\/ $\gamma,\gamma_0,\gamma_1,\gamma_2$ such that\/ $\gamma_1,\gamma_2$ are
tangent, $\alpha(\gamma_0,\gamma_1)=a_0$, $\alpha(\gamma_1,\gamma)=a_1$, $\alpha(\gamma,\gamma_2)=a_2$, and\/
$\alpha(\gamma_2,\gamma_0)=a_3$.}

\medskip

{\bf Proof.} Pick a couple of distinct tangent geodesics $\gamma_1,\gamma_2$ such that, in terms of Remark 2.9, the~triangle
$(v_1,v,v_2)$ is clockwise oriented if $a_0+a_3<\pi$ and counterclockwise oriented if $\pi<a_0+a_3$. By~Remark 2.9, there exists
a unique geodesic $\gamma_0$ with $\alpha(\gamma_2,\gamma_0)=a_3$ and $\alpha(\gamma_0,\gamma_1)=a_0$. Again by Remark 2.9, there
exists a unique geodesic $\gamma$ with $\alpha(\gamma_1,\gamma)=a_1$ and $\alpha(\gamma,\gamma_2)=a_2$ because $a_0+a_3<\pi$
implies $\pi<a_1+a_2$ and $\pi<a_0+a_3$ implies $a_1+a_2<\pi$.

\medskip

{\bf2.11.~Remark.} {\sl Let\/ $\Sigma$ be a hyperbolic\/ $2$-sphere with cone points {\rm(}$\sigma=-1${\rm)} and let\/
$R\subset\Sigma$ be a connected and simply connected region containing no cone points. Then there is a locally isometric immersion
from\/ $R$ into\/ $\B V$. If\/ $R$ is a star-like region, then the immersion is an embedding.}

\bigskip

\centerline{\bf3.~Technical lemmas}

\medskip

{\bf3.1.~Lemma.} {\sl Let\/ $k_0,k_1,k_2,k_3\in\Bbb C$ with\/ $|k_j|=1$ be fixed. Pick some\/ $c_j\in\B V$,
$0\leqslant j\leqslant3$, and denote\/ $t_1:=\ta(c_0,c_1)$, $t_2:=\ta(R_{c_1,k_1}c_0,c_2)$, and\/ $t_3:=\ta(c_0,c_3)$. If the
relation\/ $R_{c_3,k_3}R_{c_2,k_2}R_{c_1,k_1}R_{c_0,k_0}=1$ holds in\/ $\SU V$, then
$$\det\left[\smallmatrix1&t_3(k_3-\overline k_3)+\overline k_3&\big(t_2(\overline k_2-k_2)+k_2\big)\overline k_0\\t_3(\overline
k_3-k_3)+k_3&1&t_1(k_1-\overline k_1)+\overline k_1\\\big(t_2(k_2-\overline k_2)+\overline k_2\big)k_0&t_1(\overline
k_1-k_1)+k_1&1\endsmallmatrix\right]=0.\leqno{\bold{(3.2)}}$$
Conversely, given\/ $0\leqslant t_1,t_2,t_3\leqslant1$ {\rm(}the case\/ $\sigma=1${\rm)} or\/ $1\leqslant t_1,t_2,t_3$ {\rm(}the
case\/ $\sigma=-1${\rm)} satisfying\/~{\rm(3.2)}, there exist\/ $c_j$'s subject to the relations\/
$R_{c_3,k_3}R_{c_2,k_2}R_{c_1,k_1}R_{c_0,k_0}=1$ in\/ $\SU V$, $t_1=\ta(c_0,c_1)$, $t_2=\ta(R_{c_1,k_1}c_0,c_2)$, and\/
$t_3=\ta(c_0,c_3)$. Such\/ $c_0,c_1,c_2,c_3\in\B V$ are unique up to the action of\/ $\U V$ if\/ $k_j\ne\pm1$ and $t_j\ne0$ for
all\/ $j$.}

\medskip

{\bf Proof.} We assume that $\langle c_j,c_j\rangle=1$ for all $j$. Denote $a_j:=\sqrt t_j$ and choose representatives of
$c_1,c_3$ such that $\langle c_0,c_j\rangle=a_j$ for all $j=1,3$. Let $p_2:=R_{c_3,k_3}^{-1}c_0$ and $p_1:=R_{c_1,k_1}c_0$. We
choose a representative of $c_2$ such that $\langle p_1,c_2\rangle=a_2$. By Remark 2.2, the Gram matrix of the chosen
$p_2,c_3,c_0,c_1,p_1,c_2\in V$ equals
$$G:=\left[\smallmatrix1&a_3k_3&t_3(k_3-\overline k_3)+\overline k_3&&&\\a_3\overline k_3&1&a_3&&&\\t_3(\overline
k_3-k_3)+k_3&a_3&1&a_1&t_1(k_1-\overline k_1)+\overline k_1&\\&&a_1&1&a_1k_1&\\&&t_1(\overline k_1-k_1)+k_1&a_1\overline
k_1&1&a_2\\&&&&a_2&1\endsmallmatrix\right].$$
Since $R_{c_3,k_3}R_{c_2,k_2}R_{c_1,k_1}R_{c_0,k_0}c_0=c_0$ is equivalent to $R_{c_2,k_2}p_1=k_0p_2$ in these terms, we conclude
by Remark 2.2 that
$$G=\left[\smallmatrix1&a_3k_3&t_3(k_3-\overline k_3)+\overline k_3&*&\big(t_2(\overline k_2-k_2)+k_2\big)\overline
k_0&a_2\overline k_0\overline k_2\\a_3\overline k_3&1&a_3&*&*&*\\t_3(\overline k_3-k_3)+k_3&a_3&1&a_1&t_1(k_1-\overline
k_1)+\overline k_1&*\\*&*&a_1&1&a_1k_1&*\\\big(t_2(k_2-\overline k_2)+\overline k_2\big)k_0&*&t_1(\overline
k_1-k_1)+k_1&a_1\overline k_1&1&a_2\\a_2k_0k_2&*&*&*&a_2&1\endsmallmatrix\right].\leqno{\bold{(3.3)}}$$
The matrix
$M:=\left[\smallmatrix1&t_3(k_3-\overline k_3)+\overline k_3&\big(t_2(\overline k_2-k_2)+k_2\big)\overline k_0\\t_3(\overline
k_3-k_3)+k_3&1&t_1(k_1-\overline k_1)+\overline k_1\\\big(t_2(k_2-\overline k_2)+\overline k_2\big)k_0&t_1(\overline
k_1-k_1)+k_1&1\endsmallmatrix\right]$
is a submatrix of $G$. It is the Gram matrix of $p_2,c_0,p_1\in V$; since $\dim_\Bbb CV=2$, these points must be linearly
dependent and their Gram matrix, degenerate. Thus, we arrive at (3.2).

In order to prove the uniqueness, we can assume without loss of generality that the $c_j$'s are normalized so that (3.3) is the
Gram matrix of $p_2,c_3,c_0,c_1,p_1,c_2$. Note that the points $p_2,c_0,p_1\in V$ whose Gram matrix equals $M$ are unique up to
the action of $\U V$. Indeed, the claim is easy if $t_j\ne1$ for some $1\leqslant j\leqslant3$ because, in this case,
$|u_j|\ne1$, where $u_j:=t_j(\overline k_j-k_j)+k_j$. If $t_1=t_2=t_3=1$, then the points $p_2,c_0,p_1$ coincide in $\B V$ and,
again, the claim follows.

Consider the Gram matrices
$$M_3:=\left[\smallmatrix1&a_3k_3&t_3(k_3-\overline k_3)+\overline k_3\\a_3\overline k_3&1&a_3\\t_3(\overline
k_3-k_3)+k_3&a_3&1\endsmallmatrix\right],\qquad M_1:=\left[\smallmatrix1&a_1&t_1(k_1-\overline k_1)+\overline
k_1\\a_1&1&a_1k_1\\t_1(\overline k_1-k_1)+k_1&a_1\overline k_1&1\endsmallmatrix\right],$$
$$M_2:=\left[\smallmatrix1&\big(t_2(\overline k_2-k_2)+k_2\big)\overline k_0&a_2\overline k_0\overline
k_2\\\big(t_2(k_2-\overline k_2)+\overline k_2\big)k_0&1&a_2\\a_2k_0k_2&a_2&1\endsmallmatrix\right]$$
of the triples $p_2,c_3,c_0$, \ $c_0,c_1,p_1$, and $p_2,p_1,c_2$. Given $p_2,c_0,p_1\in V$ with the Gram matrix $M$, the~point
$c_j\in V$ is uniquely determined by $M_j$ if $|u_j|\ne1$. Yet, $c_j$ is uniquely determined by $M_j$ if $a_j=1$. It~remains to
observe that $|u_j|=1$ is equivalent to $t_j=0,1$ or $k_j=\pm1$.

Finally, let us show the existence of $c_j$'s. Suppose that (3.2) is valid, $\det M=0$. Since
$|u_j|^2-1=2t_j(t_j-1)(1-\Re k_j^2)$, we have $|u_j|\geqslant1$ if $\sigma=-1$ and $|u_j|\leqslant1$ if $\sigma=1$. This means
that there exist points $p_2,c_0,p_1\in V$ whose Gram matrix equals $M$. These points generate $V$ if $\rk M=2$ and coincide in
$\B V$ if $\rk M=1$. For a similar reason, the equalities $\det M_3=\det M_1=\det M_2=0$ imply that there exist $c_3,c_1,c_2$
such that $M_3,M_1,M_2$ are the Gram matrices of the triples $p_2,c_3,c_0$, \ $c_0,c_1,p_1$, and $p_2,p_1,c_2$, respectively.
From $\left[\smallmatrix k_3&a_3(1-k_3^2)&-1\endsmallmatrix\right]M_3=0$, we obtain $k_3p_2+a_3(1-k_3^2)c_3-c_0=0$, i.e.,
$R_{c_3,k_3}p_2=c_0$. From $\left[\smallmatrix k_1&a_1(\overline k_1-k_1)&-1\endsmallmatrix\right]M_1=0$, we deduce
$k_1c_0+a_1(\overline k_1-k_1)c_1-p_1=0$, i.e., $R_{c_1,k_1}c_0=p_1$. It follows from
$\left[\smallmatrix k_0&-k_2&a_2(k_2-\overline k_2)\endsmallmatrix\right]M_2=0$ that $k_0p_2-k_2p_1+a_2(k_2-\overline k_2)c_2=0$,
i.e., $R_{c_2,k_2}p_1=k_0p_2$. As $R_{c_0,k_0}c_0=\overline k_0c_0$, we obtain
$R_{c_3,k_3}R_{c_2,k_2}R_{c_1,k_1}R_{c_0,k_0}c_0=c_0$. Taking into account $R_{c_3,k_3}R_{c_2,k_2}R_{c_1,k_1}R_{c_0,k_0}\in\SU V$
and $\langle c_0,c_0\rangle\ne0$, we arrive at $R_{c_3,k_3}R_{c_2,k_2}R_{c_1,k_1}R_{c_0,k_0}=1$
$_\blacksquare$

\medskip

The relation $R_{c_3,k_3}R_{c_2,k_2}R_{c_1,k_1}R_{c_0,k_0}=1$ in $\PU V$ means
$R_{c_3,k_3}R_{c_2,k_2}R_{c_1,k_1}R_{c_0,k_0}=\pm1$ at the level of $\SU V$. Given some unitary complex numbers
$k_0,k_1,k_2,k_3\in\Bbb C$, the next two lemmas concern the solutions, at the level of $\SU V$, of the equations
$R_{c_3,k_3}R_{c_2,k_2}R_{c_1,k_1}R_{c_0,k_0}=\pm1$. (Note that, by Lemma~3.1 and Remark 2.2, equation (3.5) says that
$R_{c_3,k_3}R_{c_2,k_2}R_{c_1,k_1}R_{c_0,k_0}=1$ for some $c_0,c_1,c_2,c_3\in\B V$ while equation (3.7) says that
$R_{c_3,k_3}R_{c_2,k_2}R_{c_1,k_1}R_{c_0,k_0}=-1$ for some $c_0,c_1,c_2,c_3\in\B V$.)

\medskip

{\bf3.4.~Lemma.} {\sl Given\/ $a_0,a_1,a_2,a_3\in]0,\pi[$, we put\/ $k_j:=e^{a_ji}$ and\/ $u_j:=t_j(\overline k_j-k_j)+k_j$
for all\/ $0\leqslant j\leqslant 3$. Then the equation
$$1+2\Re(\overline k_0u_1u_2u_3)=|u_1|^2+|u_2|^2+|u_3|^2\leqno{\bold{(3.5)}}$$
in\/ $t_1,t_2,t_3\geqslant1$ has no solution with\/ $t_3=1$ if\/ $\sum_ja_j\leqslant\pi$ or\/ $3\pi\leqslant\sum_ja_j$.}

\medskip

{\bf Proof.} Since $u_3=\overline k_3$ when $t_3=1$, the equation takes the form
$2\Re(\overline k_0\overline k_3u_1u_2)=|u_1|^2+|u_2|^2$, which is equivalent to $|u_2-\overline u_1k_0k_3|^2=0$, i.e.,
to $u_2=\overline u_1k_0k_3$. Note that, replacing $a_j$ by $\pi-a_j$, we~change $k_j$ by $-\overline k_j$ and the last
equation becomes $-\overline u_2=-u_1\overline k_0\overline k_3$, i.e., it remains the same. Therefore, we may assume that
$\sum_ja_j\leqslant\pi$.

Denote $s:=\sin(a_0+a_3)$, $b:=\cos(a_0+a_3)$, $s_j:=\sin a_j$, and $b_j:=\cos a_j$ for $j=1,2$. It follows from
$\sum_ja_j\leqslant\pi$ that $s,s_1,s_2>0$. In these terms, the equation $u_2=\overline u_1k_0k_3$ takes the form
$-2s_2t_2i+b_2+s_2i=(2s_1t_1i+b_1-s_1i)(b+si)$. Hence, $b_2=b_1b+s_1s(1-2t_1)$ and $s_2(1-2t_2)=b_1s+s_1b(2t_1-1)$. This
means that $2t_1-1=\frac{b_1b-b_2}{s_1s}$ and $2t_2-1=\frac{b_2b-b_1}{s_2s}$ because $b^2+s^2=1$. Since $s_1s,s_2s>0$, from
$t_1,t_2\geqslant1$, we~conclude that $b_1b-s_1s\geqslant b_2$ and $b_2b-s_2s\geqslant b_1$, i.e.,
$\Re(k_0k_1k_3)\geqslant\Re k_2$ and $\Re(k_0k_2k_3)\geqslant\Re k_1$. In view of $\sum_ja_j\leqslant\pi$, we obtain
$a_0+a_1+a_3\leqslant a_2$ and $a_0+a_2+a_3\leqslant a_1$, implying $a_1<a_2$ and $a_2<a_1$, a contradiction
$_\blacksquare$

\medskip

{\bf3.6.~Lemma.} {\sl Given\/ $a_0,a_1,a_2,a_3\in]0,\pi[$, we put\/ $k_j:=e^{a_ji}$ and\/ $u_j:=t_j(\overline k_j-k_j)+k_j$
for all\/ $0\leqslant j\leqslant 3$. Then the solutions of the equation
$$1-2\Re(\overline k_0u_1u_2u_3)=|u_1|^2+|u_2|^2+|u_3|^2\leqno{\bold{(3.7)}}$$
in\/ $t_1,t_2,t_3\geqslant1$ constitute a compact.}

\medskip

{\bf Proof.} We take $\sigma=-1$. Let $t$ be the tance that corresponds to the distance $d$. It follows from
$\cosh^2\frac{\dist(p_1,p_2)}2=\ta(p_1,p_2)$ that $(2t-1)^2$ is the tance that corresponds to the distance $2d$.
Consequently, $64t(t-\frac12)^2(t-1)+1=\big(2(2t-1)^2-1\big)^2$ is the tance that corresponds to the distance $4d$.

By Lemma 3.1 and Remark 2.2, the equation (3.7) says that, for some points $c_0,c_1,c_2,c_3\in\B V$, the relation
$R_{c_3,k_3}R_{c_2,k_2}R_{c_1,k_1}R_{c_0,-k_0}=1$ holds in $\SU V$, where $t_1:=\ta(c_0,c_1)=\ta(c_1,p_1)$,
$t_2:=\ta(p_1,c_2)=\ta(c_2,p_2)$, $t_3:=\ta(p_2,c_3)=\ta(c_3,c_0)$, $p_1:=R_{c_1,k_1}c_0$, and $p_2:=R_{c_3,k_3}^{-1}c_0$.
Denoting the corresponding distances by $d_1,d_2,d_3$, we obtain $2d_1+2d_2\geqslant\dist(p_2,c_0)$. So, if
$d_2=\max(d_1,d_2)$, then $4d_2\geqslant\dist(p_2,c_0)$. By Remark 2.2,
$\ta(p_2,c_0)=\big|t_3(\overline k_3-k_3)+k_3\big|^2=2t_3(t_3-1)(1-\Re k_3^2)+1$. Therefore,
$32(t_2-\frac12)^2t_2(t_2-1)\geqslant t_3(t_3-1)(1-\Re k_3^2)$. Since $(t_2-\frac12)^2\geqslant t_2(t_2-1)$ and
$t_3(t_3-\nomathbreak1)\geqslant(t_3-1)^2$, we~obtain
$t_2-\frac12\geqslant(t_3-1)^{\frac12}\big(\frac{1-\Re k_3^2}{32}\big)^{\frac14}$, implying
$\max(t_1,t_2)\geqslant m(t_3-1)^{\frac12}$, where
$m:=\min\Big(\big(\frac{1-\Re k_1^2}{32}\big)^{\frac14},\big(\frac{1-\Re k_2^2}{32}\big)^{\frac14},\big(\frac{1-\Re
k_3^2}{32}\big)^{\frac14}\Big)>0$.

Without loss of generality, we may assume that $t_3\geqslant t_2\geqslant t_1\geqslant1$.

Suppose that the set of solutions of (3.7) is not compact. Then there exist solutions with arbitrary big $t_3$, i.e.,
$t_3\gg0$. The inequality $\max(t_1,t_2)\geqslant m(t_3-1)^{\frac12}$ implies $t_2\gg0$.

Denote $s_j:=\sin a_j$, $b_j:=\cos a_j$, and $x_j:=s_j(2t_j-1)$. Then $\overline k_0=b_0-s_0i$, $u_j=b_j-x_ji$, and
$x_j\geqslant s_j>0$ because $t_j\geqslant1$ and $a_j\in]0,\pi[$. Note that $t_2\gg0$ and $t_3\gg0$ imply $x_2\gg0$ and
$x_3\gg0$. In~the introduced terms, the equation (3.7) takes the form
$$1-2b_0b_1b_2b_3+2s_0(b_2b_3x_1+b_3b_1x_2+b_1b_2x_3)+2b_0(b_3x_1x_2+b_1x_2x_3+b_2x_3x_1)-2s_0x_1x_2x_3=$$
$$=b_1^2+x_1^2+b_2^2+x_2^2+b_3^2+x_3^2,$$
which can be written as
$$x_1^2+2p(x_2,x_3)x_1+q(x_2,x_3)=0,\leqno{\bold{(3.8)}}$$
where
$$p(x_2,x_3):=s_0x_2x_3-b_0b_3x_2-b_0b_2x_3-s_0b_2b_3=s_0(x_2-b_0b_2s_0^{-1})(x_3-b_0b_3s_0^{-1})-b_2b_3s_0^{-1},$$
$$q(x_2,x_3):=x_2^2-2b_0b_1x_2x_3+x_3^2-2s_0b_1b_3x_2-2s_0b_1b_2x_3+2b_0b_1b_2b_3+b_1^2+b_2^2+b_3^2-1.$$
It follows from $x_2\gg0$ and $x_3\gg0$ that $p(x_2,x_3)\gg0$. Hence, the fact that the equation (3.8) possesses a root
$x_1\geqslant s_1$ implies $s_1^2+2p(x_2,x_3)s_1+q(x_2,x_3)\leqslant0$. Consequently, we obtain
$$s_1^2+2s_0s_1x_2x_3-2b_0s_1b_3x_2-2b_0s_1b_2x_3-2s_0s_1b_2b_3+$$
$$+x_2^2-2b_0b_1x_2x_3+x_3^2-2s_0b_1b_3x_2-2s_0b_1b_2x_3+2b_0b_1b_2b_3+b_1^2+b_2^2+b_3^2-1\leqslant0,$$
which is equivalent to
$$2(1+s_0s_1-b_0b_1)x_2x_3+(x_2-x_3)^2-2(b_0s_1+s_0b_1)(b_3x_2+b_2x_3)+2(b_0b_1-s_0s_1)b_2b_3+b_1^2+b_2^2+b_3^2+s_1^2-1
\leqslant0,$$
i.e., to
$$2\big(1-\Re(k_0k_1)\big)x_2x_3+(x_2-x_3)^2-2\Im(k_0k_1)(b_3x_2+b_2x_3)+2\Re(k_0k_1)b_2b_3+b_2^2+b_3^2\leqslant0.$$
Since $1-\Re(k_0k_1)>0$ and $x_2,x_3\gg0$, we arrive at a contradiction
$_\blacksquare$

\bigskip

\centerline{\bf4.~The space $S(a_0,a_1,a_2,a_3)$ of hexagons, $\sigma=-1$}

\medskip

In this section, we study in detail the space $R(a_0,a_1,\dots,a_n)$ of relations between elliptic isometries of fixed conjugacy
classes for $n=3$.

\smallskip

From now on, we assume $\sigma=-1$ and the isometries are considered as elements in $\PU V$, unless the contrary is stated.

\medskip

{\bf4.1.~Definition.} Let $a_0,a_1,\dots,a_n\in]0,\pi[$ be given. Denote by $R:=R(a_0,a_1,\dots,a_n)$ the {\it space of
relations\/} $R_n\dots R_1R_0=1$ in $\PU V$ considered up to conjugation in $\PU V$, where the conjugacy classes of
$R_0,R_1,\dots,R_n\in\PU V$ are exactly those of $R_{c,k_0},R_{c,k_1},\dots,R_{c,k_n}$, listed perhaps in different order,
and~$k_j:=e^{a_ji}$. Even if some $a_j$'s coincide, we still consider the corresponding classes as different; in other words, the
classes are labeled with the $j$'s. Without loss of generality, we assume that the conjugacy class of $R_0$ is that of
$R_{c,k_0}$.

One can view $R$ as a {\it relative character variety,} i.e., formed by all $\PU V$-representations $\varrho:F_n\to\PU V$ of the
free group $F_n:=\langle r_0,r_1,\dots,r_n\mid r_n\dots r_1r_0=1\rangle$ of rank $n$ such that the conjugacy class of
$\varrho r_0$ is that of $R_{c,k_0}$ and the conjugacy class of $\varrho r_j$ is that of
$R_{c,k_{\beta j}}$ for all $1\leqslant j\leqslant n$, where $\beta$ is a permutation on $\{1,2,\dots,n\}$; the representations
are considered up to conjugation in $\PU V$.

Also, we can interpret $R$ as the space of labeled $2n$-gons $P$ as follows. Let $\varrho:F_n\to\PU V$ be a representation as
above and let $c_j$ stand for the fixed point of $R_j:=\varrho r_j$. The consecutive vertices of the closed piecewise geodesic
path $P$ are $c_1,p_1,c_2,p_2,\dots,c_n,p_n:=c_0$, where $p_j:=R_jp_{j-1}$ (the indices are modulo $n$). Each vertex $c_j$ is
labeled with some conjugacy class (of $R_{c,k_l}$, $l\ne0$) and $p_n$ is labeled with the conjugacy class of $R_{c,k_0}$.

\medskip

{\bf4.2.~Definition.} Denote by $m_0:F_n\to F_n$ the automorphism of $F_n$ given by $r_0\mapsto r_1r_0r_1^{-1}$,
$r_j\mapsto r_{j+1}$ for any $0<j<n$, and $r_n\mapsto r_1$; for any $0<l<n$, denote by $m_l:F_n\to F_n$ the automorphism of $F_n$
given by $r_l\mapsto r_{l+1}$, $r_{l+1}\mapsto r_{l+1}r_lr_{l+1}^{-1}$, and $r_j\mapsto r_j$ for any $j\ne l,l+1$; and denote by
$m_n:F_n\to F_n$ the automorphism of $F_n$ given by $r_0\mapsto r_1r_n^{-1}r_0r_nr_1^{-1}$, $r_1\mapsto r_1r_nr_1^{-1}$,
$r_j\mapsto r_j$ for all $1<j\leqslant n$, and $r_n\mapsto r_1$. The group $\overline M_n$ generated by $m_0,m_1,\dots,m_n$ acts
from the right on $R(a_0,a_1,\dots,a_n)$ by composition at the level of representations.

In terms of labeled $2n$-gons, the action of $m_0$ is just replacing the label $c_0$ from $p_n$ to $p_1$.

Note that $m_0m_l=m_{l+1}m_0$ for all $0<l<n$. Also, by induction on $l$, the automorphism $m_0^l$ shifts by $l$ (modulo $n$) the
indices of the $r_j$'s for $1\leqslant j\leqslant n$ and $m_0^l:r_0\mapsto(r_l\dots r_2r_1)r_0(r_l\dots r_2r_1)^{-1}$.
In~particular, $m_0^n=1$.

\medskip

Our main interest in this section is to understand when the action of $\overline M_n$ on a component of $R$ is discrete.

\medskip

{\bf4.3.~Definition.} Denote by $S:=S(a_0,a_1,\dots,a_n)$ the quotient of $R(a_0,a_1,\dots,a_n)$ by the action of the cyclic
group of order $n$ generated by $m_0$. In other words, $S$ is formed by $2n$-gons with a forgotten label $c_0$ at the $p_j$'s.

Let $1\leqslant j,l\leqslant n$ be such that $j\ne l-1,l$ (the indices are modulo $n$) and let $P\in S$ be $2n$-gon with
consecutive vertices $c_1,p_1,c_2,p_2,\dots,c_n,p_n$. The {\it modification\/} $Pm_{j,l}$ of $P$ is defined as follows.
The~geode\-sic segment $[p_j,c_l]$ cuts $P$ into two polygons $P_1$ and $P_2$ whose consecutive vertices are respectively
$p_j,c_{j+1},p_{j+1},\dots,p_{l-1},c_l$ and $c_l,p_l,c_{l+1},p_{l+1},\dots,c_j,p_j$. We rotate $P_1$ with $R_l$ so that $R_lP_1$
and $P_2$ become glued along $[c_l,p_l]$ (and remove this side), thus providing a new $2n$-gon $Pm_{j,l}\in S$ that has two new
(consecutive) sides $[p_k,c_l]$ and $R_l[c_l,p_j]=[c_l,R_lp_j]$.

The group $M_n$ generated by all modifications $m_{j,l}$ acts from the right on $S$. It is easy to see that the map $R\to S$ is a
bijection at the level of the $\overline M_n$-orbits and the $M_n$-orbits. (For a formal proof, one can use the relations
$m_0m_l=m_{l+1}m_0$.)

\medskip

Let $P\in S$ be a $2n$-gon with consecutive vertices $c_1,p_1,c_2,p_2,\dots,c_n,p_n$. Then the sides $[p_{j-1},c_j]$ and
$[c_j,p_j]$ have equal length and the `exterior' angle at $c_j$, i.e., the angle from $[c_j,p_{j-1}]$ to $[c_j,p_j]$ at $c_j$
counted in the counterclockwise sense, equals $2a_j$ (unless $p_{j-1}=c_j=p_j$). It follows from the relation $R_n\dots R_1R_0=1$
that the `interior' angles at the $p_j$'s sum to $2\pi-2a_0$ modulo $2\pi$. Conversely, if a $2n$-gon $P$ satisfies the listed
conditions, we obtain $P\in S$.

If $a_j\geqslant\frac\pi2$ for all $0\leqslant j\leqslant n$ and the `interior' angles at the $p_j$'s sum to $2\pi-2a_0$, then
$P$ is convex. In~this case, the $2n$-gon $P$ remains convex after any modification $m\in M_n$.

\medskip

{\bf4.4.~Small $n$.} Clearly, $S(a_0)=\varnothing$ because $a_0\in]0,\pi[$. By Remark 2.2, the space $S(a_0,a_1)$ consists of a
single point iff $a_0+a_1=\pi$; otherwise, $S(a_0,a_1)=\varnothing$.

\medskip

{\bf4.4.1.~Lemma {\rm(folklore)}.} {\sl Let\/ $a_0,a_1,a_2\in]0,\pi[$ be fixed. Then\/ $S(a_0,a_1,a_2)$ is nonempty exactly in
the following cases\/{\rm:}

\smallskip

$\bullet$ $\sum_ja_j=\pi$ or\/ $\sum_ja_j=2\pi$,

$\bullet$ $\sum_ja_j<\pi$,

$\bullet$ $2\pi<\sum_ja_j$.

\smallskip

In these cases, $S(a_0,a_1,a_2)$ consists of a single point and the corresponding quadrangle\/ $(c_1,p_1,c_2,p_2)$ is
respectively degenerate\/ {\rm(}i.e., $c_1=p_1=c_2=p_2${\rm)}, clockwise oriented, and counterclockwise oriented.}

\medskip

{\bf Proof.} If $c_0=c_1=c_2$, then $k_0^2k_1^2k_2^2=1$ by Remark 2.2 and $\sum_ja_j\equiv0\mod\pi$. A bit later, we will see
that the converse is also true.

It cannot happen that only two of $c_0,c_1,c_2$ coincide. So, we assume these points pairwise distinct. Denote by $I_j$ the
reflection in the geodesic $\gamma_j$ joining $c_{j-1}$ and $c_j$ (the indices are modulo $3$). Then $R_j=I_{j+1}I_j$. Indeed, it
suffices to show this for $j=1$. In the quadrangle $P=(c_1,p_1,c_2,p_2)$, the~triangles $(c_1,p_2,c_2)$ and $(c_1,p_1,c_2)$ are
congruent by means of $I_2$ because $p_1=R_1p_2\ne p_2$ in view of $a_1\in]0,\pi[$ (otherwise, $R_1$ would have two distinct
fixed points $c_1$ and $p_2=c_0$). Therefore, $I_2c_0=R_1c_0$, implying $I_2I_1c_0=R_1c_0$ and $I_2I_1c_1=R_1c_1$. The
orientation-preserving isometries $I_2I_1$ and $R_1$ coincide on two distinct points. Hence, they are equal.

By Remark 2.3, $a_j=\alpha(\gamma_j,\gamma_{j+1})$ is the oriented angle from $\gamma_j$ to $\gamma_{j+1}$. It is easy to see
that, if the triangle $(c_0,c_1,c_2)$ is clockwise (counterclockwise) oriented, then $a_j$ (respectively, $\pi-a_j$) is its
interior angle at $c_j$. So, we arrive at $\sum_ja_j<\pi$ and $2\pi<\sum_ja_j$, respectively. Obviously, under this condition,
the quadrangle $(c_1,p_1,c_2,p_2)$ is geometrically unique
$_\blacksquare$

\medskip

{\bf4.5.~Topology of $S(a_0,a_1,a_2,a_3)$.} Now, we change the convention concerning labeling the vertices $c_j$'s of a point
$P\in S(a_0,a_1,a_2,a_3)$ so that $c_j$ and $R_j$ correspond to the conjugacy class of $R_{c,k_j}$ for all
$1\leqslant j\leqslant3$. Thus, we have two types of hexagons: those with consecutive vertices $c_1,p_1,c_2,p_2,c_2,p_3$, where
$p_1=R_1p_3$, $p_2=R_2p_1$, and $p_3=R_3p_2$, and those with consecutive vertices $p_3,c_3,p_2,c_2,p_1,c_1$, where
$p_1=R_1^{-1}p_3$, $p_2=R_2^{-1}p_1$, and $p_3=R_3^{-1}p_2$. To each hexagon $P$, we associate a triple
$(t_1,t_2,t_3)\in\nomathbreak\Bbb R^3$, where $t_j:=\ta(c_j,p_j)$. In what follows, we frequently use to place the label $c_0$ at
$p_3$ thus getting the relation $R_3R_2R_1R_0=1$ for the hexagons of the first type and the relation $R_0R_1R_2R_3=1$ for the
hexagons of the second type.

According to the new convention, the group $M_3$, previously generated by $m_{1,3}$, $m_{2,1}$, $m_{3,2}$ is now generated by
the modifications $n_1$, $n_2$, $n_3$ that transform respectively the relation $R_3R_2R_1R_0=1$ into the relations
$R_0^{R_3^{-1}}R_1R_2^{R_1^{-1}}R_3=1$, $R_0^{R_3^{-1}}R_1^{R_2}R_2R_3=1$, $R_0^{R_1}R_1R_2^{R_3}R_3=1$ and the relation
$R_0R_1R_2R_3=\nomathbreak1$ into the relations $R_3R_2^{R_1}R_1R_0^{R_3}=1$, $R_3R_2R_1^{R_2^{-1}}R_0^{R_3}=1$,
$R_3R_2^{R_3^{-1}}R_1R_0^{R_1^{-1}}=1$, thus altering the type of a hexagon. It is easy to verify that $n_j^2=1$ for all
$j=1,2,3$ and that $n_j$ keeps all the $t_l$'s except perhaps the $t_j$.

We say that a point $P\in S(a_0,a_1,a_2,a_3)$ is {\it degenerate\/} with respect to $t_3$ if the only point
$P'\in S(a_0,a_1,a_2,a_3)$ of the same type as $P$ such that $t_3P'=t_3P$ is the point $P$. This implies that $t_3P=1$ or
$c_1=c_2$ for $P$ as, otherwise, we get infinitely many points with the same $t_3$ provided by $p_3,gc_1,gc_2,c_3$, where
$g\in C_R$ runs over the centralizer $C_R$ in $\PU V$ of $R:=R_0R_3$ (when $P$ is of the first type) or of $R:=R_3R_0$ (when $P$
is of the second type) and the label $c_0$ is placed at $p_3$. A point $P\in S(a_0,a_1,a_2,a_3)$ is said to be {\it boundary\/}
with respect to $t_3$ if $P$ is not degenerate with respect to $t_3$ and $t_3P=1$. A point $P\in S(a_0,a_1,a_2,a_3)$ is {\it
regular\/} with respect to $t_3$ if it is not degenerate nor boundary with respect to $t_3$.

\medskip

{\bf4.5.1.~Lemma.} {\sl A point\/ $P\in S(a_0,a_1,a_2,a_3)$ is regular with respect to\/ $t_3$ iff\/ $t_3P\ne1$ and\/
$c_1\ne c_2$ for\/ $P$.

It follows a complete list of the points in\/ $S(a_0,a_1,a_2,a_3)$ degenerate or boundary with respect to\/~$t_3${\rm:}

\smallskip

\noindent
{\bf a.}~Let\/ $\sum_ja_j<\pi$ or\/ $3\pi<\sum_ja_j$. Then there exist exactly\/ $2$ degenerate points\/ $P_1$ and\/ $P_2$ of a
given type. These points satisfy\/ $t_3P_1=1$, $c_1\ne c_2$ for\/ $P_1$, $t_3P_2\ne1$, and\/ $c_1=c_2$ for\/ $P_2$.

\noindent
{\bf b.}~Let\/ $\sum_ja_j=\pi$ or\/ $\sum_ja_j=3\pi$. Then there exists exactly\/ $1$ degenerate point of a given type. This
point satisfies $c_1=p_1=c_2=p_2=c_3=p_3$.

\noindent
{\bf c.}~Let\/ $a_0+a_3<\pi<a_1+a_2$ or\/ $a_1+a_2<\pi<a_0+a_3$. Then there exists exactly\/ $1$ degenerate point\/ $P_0$ of a
given type. The condition\/ $c_1=p_1=c_2=p_2=c_3=p_3$ for\/ $P_0$ is equivalent to\/ $\sum_ja_j=2\pi$.

\noindent
{\bf d.}~Let\/ $a_0+a_3=a_1+a_2=\pi$. Then there is no degenerate point. This case is the only one when there are points boundary
with respect to\/ $t_3$. Such points, forming a space homeomorphic to the ray\/ $[1,\infty[$, are given by the condition\/
$t_3P=1$.

\noindent
{\bf e.}~There is no degenerate point in the remaining cases, i.e., in the cases

\smallskip

$\bullet$ $\pi<\sum_ja_j$ with\/ $a_0+a_3<\pi$ and\/ $a_1+a_2<\pi$,\qquad$\bullet$ $\sum_ja_j<3\pi$ with\/ $\pi<a_0+a_3$ and\/
$\pi<a_1+a_2$,

$\bullet$ $a_0+a_3=\pi$ with\/ $a_1+a_2\ne\pi$,\hskip87pt$\bullet$ $a_1+a_2=\pi$ with\/ $a_0+a_3\ne\pi$.}

\medskip

{\bf Proof.} If $t_3P=1$ and $c_1\ne c_2$ for $P$, then, by Lemma 4.4.1, $P$ is degenerate with respect to $t_3$. If $t_3P\ne1$,
$c_1=c_2$ for $P$, and $P'$ is a point of the same type as $P$ with $t_3P'=t_3P$, then we can assume that $P$ and $P'$ share the
vertices $c_3,p_3$; so, the vertices $c_3,p_3,c_1=c_2$ determine $P$ and the vertices $c_3,p_3,c'_1,c'_2$ determine $P'$. Placing
the label $c_0$ at $p_3$, we get the relations $(R_2R_1)R_0R_3=1$ for $P$ and $R'_2R'_1R_0R_3=1$ for $P'$ in the case of the
first type, and the relations $(R_1R_2)R_3R_0=1$ for $P$ and $R'_1R'_2R_3R_0=1$ for $P'$ in the case of the second type. Since
$c_0\ne c_3$ implies $a_1+a_2\ne\pi$, the relation $R'_2R'_1(R_2R_1)^{-1}=1$ or $R'_1R'_2(R_1R_2)^{-1}=1$, respectively, provides
$c_1=c'_1=c'_2$ by Lemma 4.4.1. Thus, $P$ is degenerate with respect to $t_3$ if $t_3P=1$ and $c_1\ne c_2$ for $P$ or if
$t_3P\ne1$ and $c_1=c_2$ for $P$.

Since $t_3P\ne1$ and $c_1=c_2$ for $P$ imply $a_1+a_2\ne\pi$, by Lemma 4.4.1, there exists a (unique) point
$P\in S(a_0,a_1,a_2,a_3)$ of a given type with $t_3P\ne1$ and $c_1=c_2$ (hence, degenerate with respect to $t_3$) exactly in the
following cases

\smallskip

$\bullet$ $\sum_ja_j<\pi$,\qquad$\bullet$ $\sum_ja_j<2\pi$ and $\pi<a_1+a_2$,\qquad$\bullet$ $2\pi<\sum_ja_j$ and
$a_1+a_2<\pi$,\qquad$\bullet$ $3\pi<\sum_ja_j$.

\smallskip

Let $a_0+a_3\ne\pi$. Again by Lemma 4.4.1, the relation $R_2R_1R_0R_3=1$ or $R_1R_2R_3R_0=1$, respectively, together with
$t_3P=1$ determine uniquely a point $P\in S(a_0,a_1,a_2,a_3)$ of a given type. So, $P$ is degenerate with respect to $t_3$ when
$t_3P=1$. Moreover, by Lemma 4.4.1, there exists a (unique) point $P_0$ of a given type with $t_3P_0=1$ (hence, degenerate with
respect to $t_3$) exactly in the following cases

\smallskip

$\bullet$ $\sum_ja_j=\pi$,\qquad$\bullet$ $\sum_ja_j=2\pi$ and $a_0+a_3\ne\pi$,\qquad$\bullet$ $\sum_ja_j=3\pi$,\qquad$\bullet$
$\sum_ja_j<\pi$,

$\bullet$ $\sum_ja_j<2\pi$ and $\pi<a_0+a_3$,\qquad$\bullet$ $2\pi<\sum_ja_j$ and $a_0+a_3<\pi$,\hskip23pt$\bullet$
$3\pi<\sum_ja_j$.

\smallskip

\noindent
Note that $c_1=c_2$ with $t_3P=1$ is possible exactly in the first $3$ cases because $a_0+a_3\ne\pi$. In these cases,
$c_1=p_1=c_2=p_2=c_3=p_3$.

In the case $a_0+a_3=\pi\ne a_1+a_2$, there is no point $P$ with $t_3P=1$.

Let $a_0+a_3=a_1+a_2=\pi$. Then $t_3P=1$ implies $c_0=c_3$ and $c_1=c_2$. Now $t_1P=t_2P$ can be an arbitrary number in
$[1,\infty[$. In other words, the condition $t_3P=1$ provides a point boundary with respect to $t_3$. In this way, we listed all
possible points boundary with respect to $t_3$.

Summarizing, we arrive at the list in the lemma
$_\blacksquare$

\medskip

Assuming the label $c_0$ placed at $p_3$, to every point $P\in S(a_0,a_1,a_2,a_3)$, we associate the isometry $R:=R_0R_3$ if $P$
is of the first type and the isometry $R:=R_3R_0$, if $P$ is of the second type. It follows immediately from Lemma 4.5.1 that the
isometry $R$ is elliptic (or the identity) if $P$ is degenerate or boundary with respect to $t_3$.

In the following proposition, for any point $P\in S(a_0,a_1,a_2,a_3)$, we place the label $c_0$ at $p_3$.

\medskip

{\bf4.5.2.~Proposition.} {\sl It follows a full list of connected components of\/ $S:=S(a_0,a_1,a_2,a_3)$ of a given\/
{\rm(}either\/{\rm)} type presented with respect to the cases listed in Lemma\/ {\rm4.5.1:}

\smallskip

\noindent
{\bf a.}~In this case, there are two components. One is topologically a\/ $2$-sphere and the isometry\/ $R_0R_3$ is elliptic for
every point in this component. The other component is topologically a plane and the isometry\/ $R_0R_3$ is hyperbolic for every
point in this component.

\noindent
{\bf b.}~In this case, there are two components. One is a single point, degenerate with respect to\/ $t_3$, and the isometry\/
$R_0R_3$ is elliptic for this point. The other component is topologically a plane and the isometry\/ $R_0R_3$ is hyperbolic for
every point in this component.

\noindent
{\bf c.}~In this case, there is a unique component, topologically a plane. There are points in the component whose isometry\/
$R_0R_3$ is elliptic, parabolic, or hyperbolic.

\noindent
{\bf d.}~In this case, there is a unique component, topologically a plane. The isometry\/ $R_0R_3$ is hyperbolic or the identity
for every point in this component.

\noindent
{\bf e.}~In this case, there is a unique component, topologically a plane. The isometry\/ $R_0R_3$ is hyperbolic for every point
in this component.

\smallskip

Every component\/ $C$, except of the one consisting of a single point, contains a curve\/ $L_2\subset C$ dividing\/ $C$ into\/
$2$ parts such that\/ $t_3C=t_3L_2$ and every fibre of\/ $L_2\overset t_3\to\longrightarrow\Bbb R$ contains at most\/ $2$ points
with the unique exception in the case\/ {\bf d} where the fibre over\/ $t_3=1$ is topologically a ray. The curve\/ $L_2$ is
topologically a circle or a line when\/ $C$ is topologically a\/ $2$-sphere or a plane, respectively. The curve\/ $L_2$ admits a
smooth parameterization by\/ $t_3$ at the points where the isometry\/ $R_0R_3$ is hyperbolic.}

\medskip

{\bf Proof.} Let $C$ be a connected component of $S$. The hexagons in $C$ are all of a same type. We will deal with the first
type indicating in parentheses what happens to the case of the second one. Denote $t_0:=t_0(a_0,a_3)$ (see Remark 2.4). By Lemma
4.5.1, $t_3P<t_0$ if $P\in S$ is degenerate with respect to $t_3$.

Let $P\in C$. If $P$ is degenerate with respect to $t_3$, then the fibre of $C\overset t_3\to\longrightarrow\Bbb R$ at $P$
consists of a single point.

If $P$ is boundary with respect to $t_3$, then the fibre of $C\overset t_3\to\longrightarrow\Bbb R$ at $P$ is homeomorphic to the
ray $[1,\infty[$ by Lemma 4.5.1.

Suppose that $P$ is regular with respect to $t_3$. Denote by $\gamma_0$ the geodesic that joins $c_0$ with~$c_3$. By~Remark 2.3,
there exist unique geodesics $\gamma_1,\gamma_2$ such that $R_0=I_1I_0$ and $R_3=I_0I_2$ (such that $R_0=I_0I_1$ and $R_3=I_2I_0$
for the second type) because $R_0\ne1$ and $R_3\ne1$, where $I_j$ stands for the reflection in $\gamma_j$. Note that, by Remark
2.3, $\alpha(\gamma_0,\gamma_1)=a_0$ and $\alpha(\gamma_2,\gamma_0)=a_3$ (for the second type, $\alpha(\gamma_1,\gamma_0)=a_0$
and $\alpha(\gamma_0,\gamma_2)=a_3$).

Since $c_1$ or $c_2$ is not a fixed point of $R:=R_0R_3=I_1I_2\ne1$ (of $R:=R_3R_0=I_2I_1\ne1$ for the second type), the points
$c_0,gc_1,gc_2,c_3$, where $g$ runs over $C_R$, provide the hexagons forming the fibre of $C\overset t_3\to\longrightarrow\Bbb R$
and this fibre is homeomorphic to the centralizer $C_R$. Indeed, let $I'$ stand for the reflection in a geodesic joining $c_1$
and $c_2$ (here admitting that the points $c_1$ e $c_2$ may coincide, which is in fact impossible). Then, by Remark 2.3,
$R_1=I'I'_1$ and $R_2=I'_2I'$ ($R_1=I'_1I'$ and $R_2=I'I'_2$ for the second type) for suitable (unique) reflections in geodesics
$I'_1$ and $I'_2$. The relation $R_2R_1R_0R_3=1$ (the relation $R_1R_2R_3R_0=1$ for the second type) implies
$I'_2I'_1=I_2I_1\ne1$ (implies $I'_1I'_2=I_1I_2\ne1$ for the second type). Now the assertion follows from Remark 2.3. Indeed,
acting by $C_R$, we obtain $I'_1=I_1$ and $I'_2=I_2$ due to Remark 2.3. Since $c_1\ne c_2$, from $gc_1=c_1$and $gc_2=c_2$ for
$g\in C_R$, we conclude $g=1$.

Denote by $L_2\subset C$ the subset of all points $P\in C$ degenerate or boundary with respect to $t_3$ and of all those points
$P\in C$ regular with respect to $t_3$ that satisfy $I'_1=I_1$ and $I'_2=I_2$. By Remark 2.3, the~fibre of
$C\overset t_3\to\longrightarrow\Bbb R$ at $P\in C$ either is a single point $P\in L_2$ (when $P$ is degenerate), or is included
in $L_2$ and homeomorphic to $[1,\infty[$ (when $P$ is boundary), or is a circle and contains exactly two points in~$L_2$ (when
$P$ is regular and $R$ is elliptic), or is a line and contains exactly one point in $L_2$ (when $R$ is parabolic or hyperbolic).

Any regular point $P\in L_2$ can be described as $4$ geodesics $\gamma,\gamma_0,\gamma_1,\gamma_2$ such that
$c_0=\gamma_1\cap\gamma_0$, $c_1=\gamma\cap\gamma_1$, $c_2=\gamma_2\cap\gamma$, $c_3=\gamma_0\cap\gamma_2$,
$\alpha(\gamma_0,\gamma_1)=a_0$, $\alpha(\gamma_1,\gamma)=a_1$, $\alpha(\gamma,\gamma_2)=a_2$, $\alpha(\gamma_2,\gamma_0)=a_3$,
$R_0=I_1I_0$, $R_1=II_1$, $R_2=I_2I$, $R_3=I_0I_2$ (for the second type, $\alpha(\gamma_1,\gamma_0)=a_0$,
$\alpha(\gamma,\gamma_1)=a_1$, $\alpha(\gamma_2,\gamma)=a_2$, $\alpha(\gamma_0,\gamma_2)=a_3$,
$R_0=I_0I_1$, $R_1=I_1I$, $R_2=II_2$, $R_3=I_2I_0$), where $I,I_0,I_1,I_2$ stand respectively for the reflections in
$\gamma,\gamma_0,\gamma_1,\gamma_2$.

Suppose that there exists a point $P\in L_2$ with hyperbolic $R$. Then, by Remark 2.5, we can continuously deform $P$ in $L_2$
making $t_3$ arbitrarily big because the tance $t_3=\ta(c_0,c_3)$ is greater or equal than that between the geodesics $\gamma_1$
and $\gamma_2$. As well, we can continuously deform $P$ diminishing $t_3$ till the moment when the geodesics $\gamma_1$ and
$\gamma_2$ are becoming tangent. In other words, $]t_0,\infty[\subset t_3C$ if $R$ can be hyperbolic for a point in $L_2$.

Suppose that there exists a point $P\in L_2$ with parabolic $R$. Then, by Remark 2.9, there is a small continuous deformation of
$P$ providing points in $L_2$ with hyperbolic $R$ as well as with elliptic $R$. Therefore, $t_0\in]b,\infty[\subset t_3C$ in this
case. In terms of Remark 2.9, if the triangle $(v_1,v,v_2)$ is clockwise (counterclockwise, for the second type) oriented, then
$a_0+a_3<\pi$ and $\pi<a_1+a_2$. If the triangle $(v_1,v,v_2)$ is counterclockwise (clockwise, for the second type) oriented,
then $\pi<a_0+a_3$ and $a_1+a_2<\pi$. Hence, a point $P\in L_2$ with parabolic $R$ can exist only in the case {\bf c} of Lemma
4.5.1.

Suppose that there exists a regular point $P\in L_2$ with elliptic $R$. Denote $a:=\alpha(\gamma_2,\gamma_1)$ (denote
$a:=\alpha(\gamma_1,\gamma_2)$ for the second type) and let $c$ stand for the intersection point of the geodesics
$\gamma_1,\gamma_2$. By~Remark 2.7, when the triangle $(c,c_0,c_3)$ is clockwise (counterclockwise, for the second type)
oriented, we have $a+\pi<a_0+a_3$, and when it is counterclockwise (clockwise, for the second type) oriented, we~have
$a_0+a_3<a$. Again by Remark 2.7, we obtain $a+a_1+a_2<\pi$ or $2\pi<a+a_1+a_2$. Once again by Remark 2.7, we~can continuously
vary $a$ within the interval $]0,\pi[$, thus getting points in $L_2$, if we keep the two nonstrict inequalities, i.e., the
inequalities in one of the following $4$ variants:

\smallskip

$\bullet$ $a+\pi\leqslant a_0+a_3$ and $a+a_1+a_2\leqslant\pi$,\qquad$\bullet$ $a+\pi\leqslant a_0+a_3$ and
$2\pi\leqslant a+a_1+a_2$,

$\bullet$ $a_0+a_3\leqslant a$ and $a+a_1+a_2\leqslant\pi$,\hskip38.5pt$\bullet$ $a_0+a_3\leqslant a$ and
$2\pi\leqslant a+a_1+a_2$.

\smallskip

\noindent
(The first two pictures below illustrate the first and the fourth variants; the other two illustrate the second and the third
variants.)

\leftskip35pt
\epsfxsize=5cm
\epsfbox{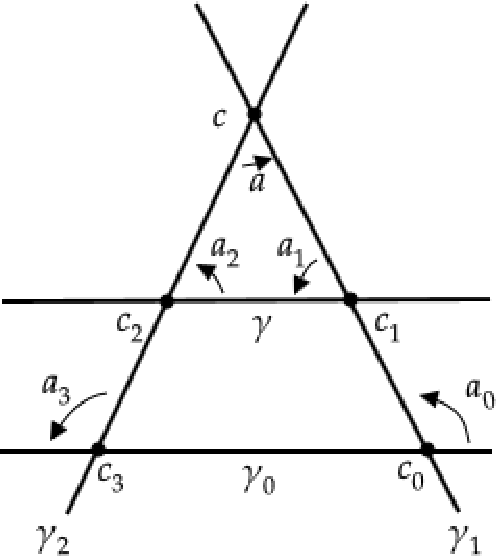}
\vskip-160pt\leftskip260pt
\epsfxsize=5cm
\noindent
\epsfbox{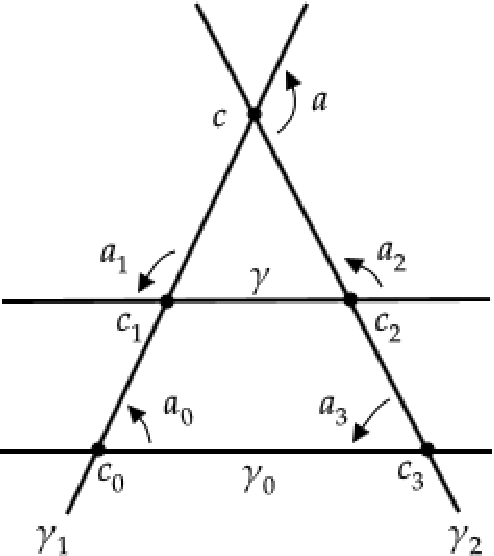}

\medskip

\leftskip35pt
\epsfxsize=5cm
\epsfbox{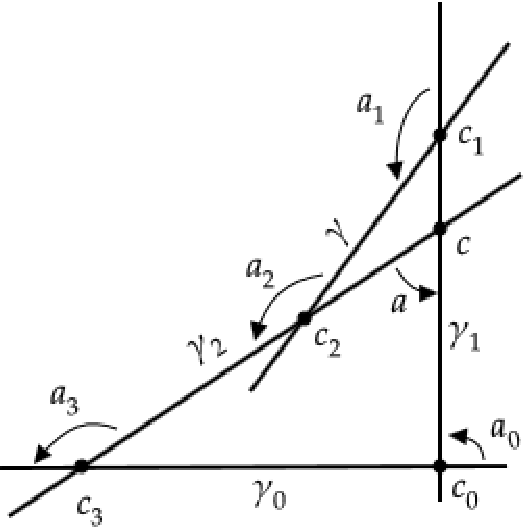}
\vskip-142pt\leftskip260pt
\epsfxsize=5cm
\noindent
\epsfbox{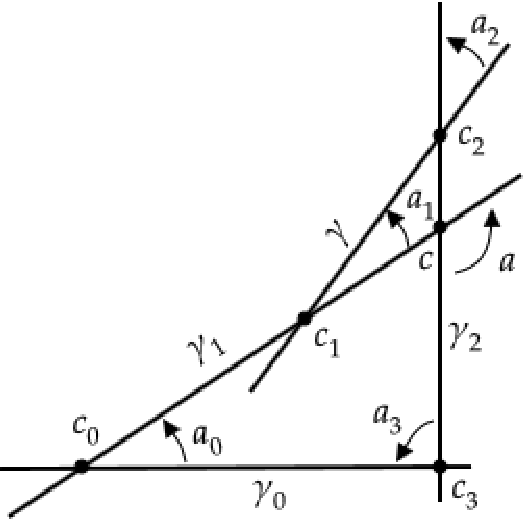}

\leftskip0pt
\bigskip
\noindent
Such a deformation can be performed by varying the tance $t_3P$ within the interval $[1,t_0[$. Note that $a$ depends
monotonically on $t_3P$ while $t_3P$ varies within the interval $]1,t_0[$ and $P\in L_2$ remains regular of the same type.
Rewriting the inequalities in a more convenient form, we arrive at

\smallskip

$\bullet$ $a\leqslant a_0+a_3-\pi$ and $a\leqslant\pi-a_1-a_2$,\qquad$\bullet$ $2\pi-a_1-a_2\leqslant a\leqslant a_0+a_3-\pi$,

$\bullet$ $a_0+a_3\leqslant a\leqslant\pi-a_1-a_2$,\hskip66.5pt$\bullet$ $a_0+a_3\leqslant a$ and $2\pi-a_1-a_2\leqslant a$.

\smallskip

\noindent
In either variant, there is some $0<a<\pi$ such that both inequalities are strict. This implies

\smallskip

$\bullet$ $a_1+a_2<\pi<a_0+a_3$ with $0<a_0+a_3-\pi<\pi$ and $0<\pi-a_1-a_2<\pi$,

$\bullet$ $3\pi<\sum_ja_j$ with $0<2\pi-a_1-a_2<a_0+a_3-\pi<\pi$,

$\bullet$ $\sum_ja_j<\pi$ with $0<a_0+a_3<\pi-a_1-a_2<\pi$,

$\bullet$ $a_0+a_3<\pi<a_1+a_2$ with $0<a_1+a_2-\pi<\pi$ and $0<\pi-a_0-a_3<\pi$.

\smallskip

\noindent
In particular, for given $a_j$'s, only one variant is possible, i.e., the four variants are disjoint.

In the second and third variants, we are in the case {\bf a} of Lemma 4.5.1 and we can reach both values
$2\pi-a_1-a_2,a_0+a_3-\pi$ and $a_0+a_3,\pi-a_1-a_2$, respectively, by varying $a$. Since the inequalities become equalities for
such points in $L_2$, we get by Remark 2.7 the two distinct degenerate points $P_1,P_2$ mentioned in Lemma 4.5.1 {\bf a}. As $C$
is connected, we obtain $[t_3P_1,t_3P_2]\subset t_3C$. We claim that $[t_3P_1,t_3P_2]=t_3C$. Indeed, otherwise, there would exist
an extra regular point $P\in L_2$ with elliptic $R$. The values of $a$ for such $P$ should fall into one of the four variants
and, since the variants are disjoint, the variant should be the same as the one we started with. This contradicts the fact that
$a$ depends monotonically on $t_3P$. Now we visualize $C$ as a topological $2$-sphere and $L_2$, as a topological circle dividing
$C$ into $2$ parts.

In the first and fourth variants, we act similarly. We can reach the values $\min(a_0+a_3-\pi,\pi-a_1-a_2)$ and
$\max(a_0+a_3,2\pi-a_1-a_2)$, respectively, by varying $a$. Since one of the inequalities becomes an equality for such a point in
$L_2$, we obtain by Remark 2.7 a degenerate point $P_0$ mentioned in Lemma~4.5.1~{\bf c}. For a similar argument,
$t_3C=[t_3P_0,b]$ or $t_3C=[t_3P_0,b[$ or $t_3C=[t_3P_0,\infty[$ for some $t_3P_0<b$. In the last case, we visualize $C$ as a
topological plane and $L_2$, as a topological line dividing $C$ into $2$ parts.

Now we may accomplish the case {\bf c}. By Remark 2.10, there exist geodesics $\gamma,\gamma_0,\gamma_1,\gamma_2$ such that
$\gamma_1,\gamma_2$ are tangent, $\alpha(\gamma_0,\gamma_1)=a_0$, $\alpha(\gamma_1,\gamma)=a_1$, $\alpha(\gamma,\gamma_2)=a_2$,
and $\alpha(\gamma_2,\gamma_0)=a_3$ (for the second type, $\alpha(\gamma_0,\gamma_2)=a_3$, $\alpha(\gamma_2,\gamma)=a_2$,
$\alpha(\gamma,\gamma_1)=a_1$, and $\alpha(\gamma_1,\gamma_0)=a_0$). Denote by $I,I_0,I_1,I_2$ the reflections in
$\gamma,\gamma_0,\gamma_1,\gamma_2$. By Remark 2.3, we obtain elliptic isometries $R_0:=I_1I_0$, $R_1:=II_1$, $R_2:=I_2I$,
$R_3:=I_0I_2$ (elliptic isometries $R_0:=I_0I_1$, $R_1:=I_1I$, $R_2:=II_2$, $R_3:=I_2I_0$ for the second type) of the conjugacy
classes $R_{c,k_0}$, $R_{c,k_1}$, $R_{c,k_2}$, $R_{c,k_3}$ subject to $R_2R_1R_0R_3=1$ (to $R_1R_2R_3R_0=1$ for the second type)
with parabolic $R$. So, we have constructed a point $P\in S$ with parabolic $R$. This point belongs to some component $C$. Hence,
there is a regular point $P'\in L_2\subset C$ with elliptic $R$. As we saw earlier, $C$ should contain the degenerate point $P_0$
considered above. Therefore, such a component (of a given type) is~unique and we are done in the case {\bf c} because we already
know that the presence of a point $P\in L_2$ with hyperbolic $R$ (induced by the presence of a point with parabolic $R$) implies
$]t_0,\infty[\subset t_3C$.

By Remarks 2.6, 2.4, and 2.3, in an arbitrary case {\bf a--e}, we get in a similar way a point $P_t\in S$ of a given type,
depending smoothly on $t>t_0$, with hyperbolic $R$ such that $t_3P_t=t$. This point generates a connected component $C$ of $S$
such that $P_t\in L_2$ and we already know that $t_3C\supset]t_0,\infty[$. Hence, by~the uniqueness in Remark 2.6, there exists a
unique connected component $C$ of $S$ (of a given type) containing a point with hyperbolic $R$. In the case {\bf c}, such a
component was already constructed above. In the case {\bf d}, we obtain $t_3C=[t_0,\infty[$ by Lemma 4.5.1. In all the other
cases, we have $t_3C=]t_0,\infty[$ because $t_0\in t_3C$ would imply the existence of a point with parabolic $R$ and we have
already seen that this is possible only in the case {\bf c}. Moreover, as $P_t$ depends smoothly on $t>t_0$, the curve $L_2$ gets
a smooth parameterization by $t_3$ at the points with hyperbolic $R$.

We have now accomplished the cases {\bf b} and {\bf e}. Indeed, in these cases, there cannot exist a regular point $P\in L_2$
with parabolic or elliptic $R$ since this is possible only in the cases {\bf c} and {\bf a}, respectively. As $t_3P'<t_0$ for any
point $P'\in S$ degenerate with respect to $t_3$, we conclude that, in the case {\bf b}, there exist a unique connected component
$C$ of $S$, $t_3C=]t_0,\infty[$, whose points have hyperbolic $R$ and a unique component consisting of a single degenerate point.
Similarly, in the case {\bf e}, we get a unique connected component $C$ of $S$, $t_3C=]t_0,\infty[$, whose points have hyperbolic
$R$. Thus, in the cases {\bf b} and {\bf e}, we~visualize the unique connected component $C$ of $S$ containing points with
hyperbolic $R$ as a topological plane and $L_2$, as a topological line dividing $C$ into $2$ parts.

There is a unique connected component $C$, $]1,\infty[\subset t_3C$, in the case {\bf d}. By Lemma 4.5.1, $t_3C=[1,\infty[$.
Clearly, $(t_1P,t_2P,t_3P)$ tends to $(t_1,t_2,1)$ with $t_1=t_2$ when a point $P\in C$ tends to a point boundary with respect to
$t_3$ because $c_0=c_3$ and $c_1=c_2$ for a boundary point. Given $t\in]1,\infty[$, there exists a unique point $P_t\in L_2$,
depending continuously on $t$, such that $t_3P_t=t$. Any point $P\in C$ with $t_3P=t$ has the form $P=gP_t$ for a suitable
(unique) isometry $g\in C_R$. Note that, for $t$ close to $1$, we have $t_1(gP_t)$ close to $t_1(g^{-1}P_t)$ just because $c_0$
and $c_1$ are close for $P_t$. This means that the points $gP_t$ and $g^{-1}P_t$ tend to a same boundary point when $gP_t$ tends
to a boundary point. In other words, in the fibre over $t_3=1$, two distinct points over $t>1$ tend to collide when $t$ tends to
$1$ (unless $g$ tends to $1$). So, in the case {\bf d}, we visualize $C$ as a topological plane and $L_2$, as a topological line
dividing $C$ into $2$ parts.

Finally, in the case {\bf a}, it suffices to construct a regular point $P\in S$ with elliptic $R$ because, as we have seen above,
such a point $P$ provides the two degenerate points $P_1$ and $P_2$ in the connected component containing $P$, implying the
existence and uniqueness of a connected component of $S$ whose points have elliptic $R$.

By Remark 2.8, there exist geodesics $\gamma,\gamma_0,\gamma_1,\gamma_2$ such that $\gamma_1,\gamma_2$ intersect,
$\alpha(\gamma_0,\gamma_1)=a_0$, $\alpha(\gamma_1,\gamma)\allowmathbreak=a_1$, $\alpha(\gamma,\gamma_2)=a_2$, and
$\alpha(\gamma_2,\gamma_0)=a_3$ (for the second type, $\alpha(\gamma_0,\gamma_2)=a_3$, $\alpha(\gamma_2,\gamma)=a_2$,
$\alpha(\gamma,\gamma_1)=a_1$, and $\alpha(\gamma_1,\gamma_0)=a_0$). Denote by $I,I_0,I_1,I_2$ the reflection in
$\gamma,\gamma_0,\gamma_1,\gamma_2$. By Remark 2.3, the isometries $R_0:=I_1I_0$, $R_1:=II_1$, $R_2:=I_2I$, $R_3:=I_0I_2$ (the
isometries $R_0:=I_0I_1$, $R_1:=I_1I$, $R_2:=II_2$, $R_3:=I_2I_0$ for the second type) of the conjugacy classes $R_{c,k_0}$,
$R_{c,k_1}$, $R_{c,k_2}$, $R_{c,k_3}$ subject to $R_2R_1R_0R_3=1$ (to $R_1R_2R_3R_0=1$ for the second type) provide a regular
point $P\in S$ with elliptic $R$
$_\blacksquare$

\medskip

For given $a_0,a_1,a_2,a_3\in]0,\pi[$, we denote $k_j:=e^{a_ji}\ne\pm1$ and $u_j:=t_j(\overline k_j-k_j)+k_j$.

Let $C\subset S:=S(a_0,a_1,a_2,a_3)$ be a connected component of $S$ and let $P\in S$. Then, placing the label $c_0$ at $p_3$,
the relation $R_{c_3,k_3}R_{c_2,k_2}R_{c_1,k_1}R_{c_0,k_0}=1$ holds in $\PU V$ if $P$ is of the first type and the relation
$R_{c_1,k_1}R_{c_2,k_2}R_{c_3,k_3}R_{c_0,k_0}=1$ holds in $\PU V$ if $P$ is of the second type. At the level of $\SU V$, we
obtain respectively the relations $R_{c_3,k_3}R_{c_2,k_2}R_{c_1,k_1}R_{c_0,k_0}=\pm1$ and
$R_{c_1,k_1}R_{c_2,k_2}R_{c_3,k_3}R_{c_0,k_0}=\pm1$ in $\SU V$. They can be rewritten as the relations
$R_{c_3,k_3}R_{c_2,k_2}R_{c_1,k_1}R_{c_0,\pm k_0}=1$ and $R_{c_1,k_1}R_{c_2,k_2}R_{c_3,k_3}R_{c_0,\pm k_0}=1$ in $\SU V$. By
Lemma 3.1, we arrive at the equation (3.5) (when $\pm k_0=k_0$) or at the equation (3.7) (when $\pm k_0=-k_0$) independently of
the type of $P$. This means that $C$ satisfies the equation (3.5) or~the equation (3.7).

The following corollary claims the converse, i.e., that the components of $S$ are given by the type of hexagons and by one of the
equations (3.5) and (3.7).

\medskip

{\bf4.5.3.~Corollary.} {\sl Let\/ $a_0,a_1,a_2,a_3\in]0,\pi[$. Denote\/ $k_j:=e^{a_ji}$ and\/ $u_j:=t_j(\overline k_j-k_j)+k_j$.

For any type of hexagons, the solutions of the equation\/ {\rm(3.5)} in\/ $t_1,t_2,t_3\geqslant1$ form a connected component\/
$C_1$ of\/ $S:=S(a_0,a_1,a_2,a_3)$ of this type, topologically a plane, with hyperbolic isometry\/ $R_0R_3$ for some point in the
component.

The equation\/ {\rm(3.7)} in\/ $t_1,t_2,t_3\geqslant1$ has no solutions in the cases\/ {\bf c--e} of Lemma\/ {\rm4.5.1.} In the
case\/~{\bf b}, the equation\/ {\rm(3.7)} in\/ $t_1,t_2,t_3\geqslant1$ has a unique solution that constitutes a connected
component\/ $C_0$ of\/ $S$ of any given type and the isometry\/ $R_0R_3$ is elliptic for this point. In the case\/ {\bf a}, the
solutions of the equation\/ {\rm(3.7)} in\/ $t_1,t_2,t_3\geqslant1$ form a connected component\/ $C_0$ of\/ $S$ {\rm(}of either
type\/{\rm)}, topologically a\/ $2$-sphere, and the isometry\/ $R_0R_3$ is elliptic for every point in this component.

The above is a full list of connected components of\/ $S$.}

\medskip

{\bf Proof.} Let us first observe that a point $P\in S$ degenerate with respect to $t_3$ does not satisfy the equation (3.5) in
the cases {\bf a--b}. Indeed, if $t_3P=1$, this assertion is just Lemma 3.4. If $c_1=c_2$, we~can treat the point $P$ as a point
in the space $S':=S(a_2,a_3,a_0,a_1)$: if $P$ corresponds to the relation $R_3R_2R_1R_0=1$ (to the relation $R_0R_1R_2R_3=1$ for
the second type), then $P'$ corresponds to the relation $R_1R_0R_3R_2=1$ (to the relation $R_2R_3R_0R_1=1$ for the second type),
and the assertion follows from Lemma 3.4 applied to $S'$, where $t_3P'$ becomes $1$, because the equation (3.5) and the cases
{\bf a--b} remain the same for $S'$.

In an arbitrary case {\bf a--e}, there exists a unique connected component $C_1$ of $S$ (of either type) possessing a point with
hyperbolic $R_0R_3$. Such a component $C_1$ cannot be compact since $]t_0,\infty[\subset t_3C_1$. By~Lemma~3.6, $C_1$ has to
satisfy the equation (3.5). By Proposition 4.5.2, there is a unique connected component (of~either type) in the cases {\bf c--e};
so, we are done in these cases.

In the remaining cases {\bf a--b}, there are two connected components $C_0,C_1$ (of either type) and $C_0$, the compact one, does
not satisfy the equation (3.5) because, by Lemma 4.5.1, it contains a point $P$ degenerate with respect to $t_3$. On the other
hand, the component $C_1$ containing a point with hyperbolic $R_0R_3$ is not compact because $]t_0,\infty[\subset t_3C_1$, hence,
$C_1$ cannot satisfy the equation (3.7) by Lemma 3.6. Consequently, $C_0$ satisfies (3.7) and $C_1$ satisfies (3.5) and we are
done
$_\blacksquare$

\medskip

We already know (see the second paragraph in 4.5) that the group $M_3$ is generated by the modifications $n_1,n_2,n_3$ and that
$n_j$ alters the type of a hexagon and keeps $t_l$ for $l\ne j$. Let us write the equations~(3.5) and (3.7) in the form
$e_1(t_1,t_2,t_3)=0$ and $e_0(t_1,t_2,t_3)=0$, respectively. Since $n_j$ transforms the relation
$R_{c_3,k_3}R_{c_2,k_2}R_{c_1,k_1}R_{c_0,k_0}=1$ in $\SU V$ into the relation $R_{c_1,k_1}R_{c_2,k_2}R_{c_3,k_3}R_{c_0,k_0}=1$ in
$\SU V$ and {\sl vice-versa\/} (analogously, for the relation $R_{c_3,k_3}R_{c_2,k_2}R_{c_1,k_1}R_{c_0,-k_0}=1$ in $\SU V$), this
modification preserves the equations $e_l(t_1,t_2,t_3)=0$ by Lemma 3.1. It is easy to see that
$e_l(t_1,t_2,t_3)=4t_j^2\sin^2a_j+\text{lower terms in }t_j\ $ with $\sin a_j\ne0$, i.e., the equation $e_l(t_1,t_2,t_3)=0$ is
quadratic in $t_j$. Therefore, at~the level of the $t_k$'s, the modification $n_j$ simply interchanges the roots of
$e_l(t_1,t_2,t_3)=0$ in $t_j$. In what follows, we may not distinguish anymore the types of hexagons and we may identify the
hexagons of different types with the same $(t_1,t_2,t_3)$. In this way, $n_j$ acts as an involution on every component $C_l$
of~$S$, where $C_l$ is now given by the equation $e_l(t_1,t_2,t_3)=0$. (Perhaps, it would be more correct to pass to the index
$2$ subgroup $M$ of the group $M_3$. However, it is more convenient to deal with the involutions and this makes no difference in
the problem of the discreteness of the action of $M_3$ on~$S$.)

\medskip

{\bf4.5.4.~Definition.} Let $1\leqslant j\leqslant 3$. Denote by $A_j\subset S:=S(a_0,a_1,a_2,a_3)$ the set of fixed points of
the involution $n_j$. By the above, in the connected component $C_l$ of $S$ given by the equation $e_l(t_1,t_2,t_3)=0$, the set
$A_j$ is given by the equation $\frac{\partial}{\partial t_j}e_l(t_1,t_2,t_3)=0$. In particular, $A_1\cap A_2\cap A_3$ is the set
of all singular points of $S$.

In the proof of Proposition 4.5.2, placing the label $c_0$ at $p_3$, we defined a curve $L_2\subset S$ of all points $P\in S$
degenerate or boundary with respect to $t_3$ and of all those points $P\in S$ regular with respect to $t_3$ that are subject to
$R_0=I_1I_0$, $R_1=II_1$, $R_2=I_2I$, $R_3=I_0I_2$ when $P$ is of the first type and are subject to $R_0=I_0I_1$, $R_1=I_1I$,
$R_2=II_2$, $R_3=I_2I_0$ when $P$ is of the second type, where $I,I_0,I_1,I_2$ are suitable reflections in geodesics.

\medskip

{\bf4.5.5.~Lemma.} {\sl$A_2=L_2$.}

\medskip

{\bf Proof.} If $P\in S$ is degenerate with respect to $t_3$, then $P$ is clearly a fixed point of $n_2$ because $P$ is a unique
point in $S$ with the value $t_3P$ of $t_3$ and $n_2$ keeps the values of $t_3$.

Definition 4.5.4 says that $P\in A_2$ iff $\ta(c_2,p_2)=\ta(p_3,c_2)$ for $P$. For a point $P\in S$ boundary with respect to
$t_3$, this condition is empty since $p_2=p_3$ in this case.

So, we assume $P\in S$ to be regular with respect to $t_3$. Then $c_0\ne c_3$ (we place the label $c_0$ at $p_3$) and
$c_1\ne c_2$.

Suppose that $P\in L_2$. Then $R_1=II_1$ if $P$ is of the first type and $R_1=I_1I$ if $P$ is of the second type, where
$Ic_2=c_2$, $p_1=R_1p_3=II_1c_0=Ic_0$ for $P$ of the first type, and $p_1=R_1^{-1}p_3=II_1c_0=Ic_0$ for $P$ of the second type.
Therefore, $\ta(c_2,p_2)=\ta(p_1,c_2)=\ta(Ic_0,Ic_2)=\ta(p_3,c_2)$, i.e., $P\in A_2$.

Conversely, let $P\in A_2$, i.e., $\ta(c_2,p_2)=\ta(p_3,c_2)$. Denote by $I$ and $I_0$ the reflections in the geodesics joining
$c_1,c_2$ and $c_0,c_3$, respectively. Since $\ta(c_0,c_1)=\ta(c_1,p_1)$ and $\ta(p_1,c_2)=\ta(c_2,p_2)=\ta(c_0,c_2)$, the
triangles $(c_1,c_0,c_2)$ and $(c_1,p_1,c_2)$ are congruent. As $p_1=c_0$ implies $c_0=c_1=p_1$ due to $p_1=R_1c_0$ (when $P$ is
of the first type) or $c_0=R_1p_1$ (when $P$ is of the second type), we obtain $p_1=Ic_0$.

Suppose that $c_2\ne c_3$ and denote by $I_2$ the reflection in the geodesic joining $c_2$ and $c_3$. Since
$\ta(c_2,p_2)=\ta(c_0,c_2)$ and $\ta(p_2,c_3)=\ta(c_3,c_0)$, the triangles $(c_2,p_2,c_3)$ and $(c_2,c_0,c_3)$ are congruent. As
$c_0=R_3p_2$ or $p_2=R_3c_0$, the equality $p_2=c_0$ would imply $p_2=c_0=c_3$. Hence, $p_2=I_2c_0$ and $p_2\ne c_2$ because
$p_2=c_2$ and $\ta(c_2,p_2)=\ta(c_0,c_2)$ would imply $p_2=c_0$.

Let $P$ be of the first type. Then $R_3p_2=c_0$, $R_3c_3=c_3$, $I_0I_2p_2=I_0I_2I_2c_0=c_0$, $I_0I_2c_3=c_3$, $R_2p_1=p_2$,
$R_2c_2=c_2$, $I_2Ip_1=I_2IIc_0=I_2c_0=p_2$, $I_2Ic_2=c_2$ with $c_0\ne c_3$ and $p_2\ne c_2$. Hence, $R_3=I_0I_2$ and
$R_2=I_2I$. Writing $R_0$ as $R_0=I_1I_0$, where $I_1$ is a suitable reflection in a geodesic, we~conclude from the relation
$R_3R_2R_1R_0=1$ that $R_1=II_1$. Therefore, $P\in L_2$.

Let $P$ be of the second type. Then $R_3c_0=p_2$, $R_3c_3=c_3$, $I_2I_0c_0=I_2c_0=p_2$, $I_2I_0c_3=c_3$, $R_2p_2=p_1$,
$R_2c_2=c_2$, $II_2p_2=II_2I_2c_0=Ic_0=p_1$, $II_2c_2=c_2$ with $c_0\ne c_3$ and $p_2\ne c_2$. Hence, $R_3=I_2I_0$ and
$R_2=II_2$. Writing $R_0$ as $R_0=I_0I_1$, where $I_1$ is a suitable reflection in a geodesic, we conclude from the relation
$R_1R_2R_3R_0=1$ that $R_1=I_1I$. Therefore, $P\in L_2$.

Suppose that $c_2=c_3$ and $c_0\ne c_1$. Denote by $I_1$ the geodesic joining $c_0$ and $c_1$.

If $P$ is of the first type, then $R_1=II_1$ because $R_1c_0=p_1$, $R_1c_1=c_1$, $II_1c_0=Ic_0=p_1$, $II_1c_1=c_1$ with
$c_0\ne c_1$. We have $R_2=I_2I$ and $R_3=I'_0I_2$, where $I_2$ and $I'_0$ are suitable reflections in geodesics. From the
relation $R_3R_2R_1R_0=1$, we obtain $R_0=I_1I'_0$, hence, $I'_0c_0=c_0$. Since $c_0\ne c_3$ and $I'_0c_3=c_3$, we conclude that
$I'_0=I_0$. Therefore, $P\in L_2$.

If $P$ is of the second type, then $R_1=I_1I$ because $R_1p_1=c_0$, $R_1c_1=c_1$, $I_1Ip_1=I_1IIc_0=c_0$, $I_1Ic_1=c_1$ with
$c_0\ne c_1$. We have $R_2=II_2$ and $R_3=I_2I'_0$, where $I_2$ and $I'_0$ are suitable reflections in geodesics. From the
relation $R_1R_2R_3R_0=1$, we obtain $R_0=I'_0I_1$, hence, $I'_0c_0=c_0$. Since $c_0\ne c_3$ and $I'_0c_3=c_3$, we conclude that
$I'_0=I_0$. Therefore, $P\in L_2$.

Finally, suppose that $c_0=c_1$ and $c_2=c_3$. Then $I=I_0$. Since $c_0\ne c_3$, we get $R_1R_0=1$ and $R_3R_2=1$. Since
$R_0=I_1I_0$, $R_3=I_0I_2$ for $P$ of the first type and $R_0=I_0I_1$, $R_3=I_2I_0$ for $P$ of the second type, where $I_1$ and
$I_2$ are suitable reflections in geodesics, we arrive at $R_1=I_0I_1$, $R_2=I_2I_0$ and at $R_1=I_1I_0$, $R_2=I_0I_2$,
respectively. Therefore, $P\in L_2$
$_\blacksquare$

\medskip

{\bf4.5.6.~Lemma.} {\sl A point\/ $P\in S:=S(a_0,a_1,a_2,a_3)$ belongs to\/ $A_1\cap A_2$ iff\/ $P$ is degenerate or boundary
with respect to\/ $t_3$.

A point\/ $P\in S$ is singular iff\/ $c_1=p_1=c_2=p_2=c_3=p_3$ for\/ $P$. Therefore, $S$ is not smooth iff\/ $\sum_ja_j=2\pi$. In
this case, there exists a unique singular point in\/ $S${\rm;} it belongs to the noncompact component\/ $C_1$.}

\medskip

{\bf Proof.} We observe first that there exists at most one point equidistant from $3$ pairwise distinct points. Clearly,
$P\in A_j$ iff $c_j$ is equidistant from $p_1,p_2,p_3$.

Let $P\in A_1\cap A_2$. If the points $p_1,p_2,p_3$ are pairwise distinct, then $c_1=c_2$. If $p_1=p_2$, then $p_1=c_2=p_2$,
hence, $p_3=c_2$, therefore, $c_1=p_1=c_2$ again. If $p_1=p_3$, then $p_1=c_1=p_3$, hence, $p_2=c_1$, therefore, $c_1=c_2$ again.
If $p_2=p_3$, then $c_3=p_3$. So, $P$ is degenerate or boundary with respect to $t_3$. It is easy to see that the converse is
also true: if $c_1=c_2$ or $c_3=p_3$ for $P$, then $P\in A_1\cap A_2$.

We know now that $P$ is singular iff $c_1=c_2$ or $c_3=p_3$ and, simultaneously, $c_2=c_3$ or $c_1=p_1$. If~$c_1=c_2=c_3$, then,
placing the label $c_0$ at $p_3$, it follows from $R_3R_2R_1R_0=1$ or from $R_1R_2R_3R_0=1$ that $c_0=c_1$ because $R_0\ne1$.
Therefore, $c_1=p_1=c_2=p_2=c_3=p_3$. If $c_1=p_1=c_2$ (or $c_2=c_3=p_3$), then $p_3=c_1$ (respectively, $p_2=c_3$) and, again,
$c_1=p_1=c_2=p_2=c_3=p_3$. If $c_1=p_1$ and $c_3=p_3$, then $p_3=c_1$ and $p_2=c_3$, implying $c_1=p_1=c_2=p_2=c_3=p_3$ because
$R_2$ has a unique fixed point
$_\blacksquare$

\medskip

{\bf4.6.~Discreteness of $M_3$.} Placing the label $c_0$ at $p_3$, an arbitrary hexagon $P\in S(a_0,a_1,a_2,a_3)$ of~the first type
corresponds to a relation $R_3R_2R_1R_0=1$. It is easy to see that, the hexagon $Pn_2n_1$ corresponds to the relation
$R_3R_2^{(R_1^{R_2})}R_1^{R_2}R_0=1$, i.e, to the relation $R_3R_2^{R_2R_1}R_1^{R_2R_1}R_0=1$. By~induction on $l$, $P(n_2n_1)^l$
corresponds to the relation $R_3R_2^{(R_2R_1)^l}R_1^{(R_2R_1)^l}R_0=1$. Since $R_0R_3=(R_2R_1)^{-1}$, the~action of $n_2n_1$
preserves the fibres of $S(a_0,a_1,a_2,a_3)\overset t_3\to\longrightarrow\Bbb R$.

In particular, if $P$ is regular with respect to $t_3$ with elliptic $R_0R_3$ and the action of $M_3$ on the connected component
containing $P$ is discrete, then $R_0R_3$ has to be periodic and the period has to be the same when we slightly vary $t_3P$. This
is impossible because the triangle $(c_0,c,c_3)$ that corresponds to the relation $R_0^{-1}(R_0R_3)R_3^{-1}=1$ (see, for
instance, the proof of Lemma 4.4.1) is completely determined by its angles and, hence, does not admit any variation of $t_3$,
where $c$ stands for the fixed point of $R_0R_3$. Consequently, the group $M$ cannot act discretely on a component containing a
point regular with respect to $t_3$ with elliptic $R_0R_3$. (We denote by $M$ the index $2$ subgroup in the group $M_3$.)

\medskip

{\bf4.6.1.~Lemma.} {\sl If\/ $M$ act discretely on a connected component of\/ $S(a_0,a_1,a_2,a_3)$, then the component is the
noncompact one\/ {\rm(}we drop the component of a single point\/{\rm)} and

\smallskip

$\bullet$ $a_j+a_k\leqslant\pi$ for all\/ $j\ne k$ or\/ $\pi\leqslant a_j+a_k$ for all\/ $j\ne k$.}

\medskip

{\bf Proof.} Suppose that $a_j+a_k<\pi<a_l+a_m$ with $j\ne k$ and $l\ne m$. Assuming without loss of generality that
$a_0\leqslant a_3\leqslant a_1\leqslant a_2$, we obtain $a_0+a_3<\pi<a_1+a_2$, i.e., we are in the case {\bf c}. By~Proposition
4.5.2, we get a unique component. It possesses a point regular with respect to $t_3$ with elliptic $R_0R_3$. As we saw above, $M$
cannot act discretely on this component
$_\blacksquare$

\medskip

In the sequel, we deal only with the cases listed in Lemma 4.6.1 and with the noncompact connected component $C_1$.

By Lemma 4.5.6, the component $C_1$ is not smooth only in the case $a_0=a_1=a_2=a_3=\frac\pi2$. Since, by~Lemma 4.5.6,
$A_1\cap A_2\cap C_1\ne\varnothing$ only in the case {\bf d} of Lemma 4.5.1, the curves $A_j$ are pairwise disjoint in $C_1$ if
$C_1$ is smooth.

For $1\leqslant j\leqslant 3$ and $P\in S(a_0,a_1,a_2,a_3)$, denote $t'_jP:=\ta(c_j,p_{j+1})$ (the indices are modulo $3$).
So,~the~condition $t_jP=t'_jP$ is equivalent to $P\in A_j$. By Proposition 4.5.2 and Lemma 4.5.5, assuming $C_1$ smooth, the
curve $A_j$ is a smooth line in $C_1$ and divides the smooth plane $C_1$ into two parts $H_j$ and $H'_j$ given by the
inequalities $t_j\leqslant t'_j$ and $t_j\geqslant t'_j$, respectively.

\medskip

{\bf4.6.2.~Lemma.} {\sl Suppose that

\smallskip

$\bullet$ $a_j+a_k\leqslant\pi$ for all\/ $j\ne k$ or\/ $\pi\leqslant a_j+a_k$ for all\/ $j\ne k$.

\smallskip

\noindent
Then the noncompact connected component\/ $C_1$ of $S(a_0,a_1,a_2,a_3)$, i.e., the one given by the equation\/
$e_1(t_1,t_2,t_3)=0$, is smooth unless\/ $a_0=a_1=a_2=a_3=\frac\pi2$. Assuming\/ $C_1$ smooth, the curves\/ $A_1,A_2,A_3$ in\/
$C_1$ are smooth pairwise disjoint lines, each divides the smooth plane $C_1$ into two parts, and, moreover,
$A_k\cap C_1\subset H_j\cap C_1$ for all\/ $j$ and\/ $k$.}

\medskip

{\bf Proof.} Suppose that, say, $A_2\cap C_1\not\subset H_1\cap C_1$. Note that the inequalities $0<a_j<\pi$ and
$a_j+a_k\leqslant\pi$ for all $j\ne k$ define in $\Bbb R^4$ a convex region that remains connected after removing the point
$(\frac\pi2,\frac\pi2,\frac\pi2,\frac\pi2)$ (similarly, for the inequalities $\pi\leqslant a_j+a_k$). Thus, we can vary the
$a_j$'s keeping $C_1$ smooth. Since $A_1\cap A_2\cap C_1=\varnothing$ during the deformation and the equation
$e_1(t_1,t_2,t_3)=0$ as well as the functions $t'_1,t'_2$ depend continuously on the $a_j$'s, we obtain
$A_2\cap C_1\not\subset H_1\cap C_1$ for the special case of $a_1=a_2=a_3=\frac\pi2$ and $a_0\in]0,\frac\pi2[$. In this case, a
hexagon are just a triangle with the vertices $p_1,p_2,p_3$ whose interior angles sum to $2a_0$, and $c_1,c_2,c_3$ are
respectively the middle points of the sides $[p_3,p_1]$, $[p_1,p_2]$, $[p_2,p_3]$ of the triangle. As
$A_2\cap C_1\not\subset H_1\cap C_1$ implies $A_2\cap H_1\cap C_1=\varnothing$, the isosceles triangle $(p_1,p_2,p_3)$ with
$\dist(p_1,p_3)=\dist(p_2,p_3)$ and $\dist(p_1,p_2)=2\dist(c_2,p_3)$ whose interior angles sum to $2a_0$ satisfies
$2\dist(c_1,p_2)<\dist(p_1,p_3)$, which is impossible
$_\blacksquare$

\medskip

{\bf4.6.3.~Theorem.} {\sl The group\/ $M_3$ generated by the modifications acts discretely on a component\/ $C$ of\/
$S(a_0,a_1,a_2,a_3)$ iff\/ $C$ is noncompact,

\smallskip

$\bullet$ $a_j+a_k\leqslant\pi$ for all\/ $j\ne k$ or\/ $\pi\leqslant a_j+a_k$ for all\/ $j\ne k$,

\smallskip

\noindent
and\/ $a_j\ne\frac\pi2$ for some\/ $0\leqslant j\leqslant3$. In this case, the quotient\/ $C/M$ is a\/ $3$-holed\/ $2$-sphere
and\/ $M$ is a free group of rank\/ $2$, where\/ $M$ stands for the index\/ $2$ subgroup in\/ $M_3$.}

\medskip

{\bf Proof.} Suppose that $a_j\ne\frac\pi2$ for some $0\leqslant j\leqslant3$. By Lemma 4.6.1, it suffices to show that $M_3$
acts discretely on the noncompact component $C_1$ in the cases mentioned in the theorem.

Denote $T:=H_1\cap H_2\cap H_3$. By Lemma 4.6.2 and Proposition 4.5.2, the copies $Tn_1,Tn_2,Tn_3$ are pairwise disjoint and
intersect $T$ only in the lines $A_1,A_2,A_3$, respectively. Moreover, $A_j$ lies in the interior of
$T_1:=T\cup Tn_1\cup Tn_2\cup Tn_3$ and $T_1$ is closed, connected, and limited by the pairwise disjoint lines
$A_3n_1,A_2n_1,A_1n_2,A_3n_2,A_2n_3,A_1n_3$. Applying to $T_1$ the involutions
$n_1n_3n_1,n_1n_2n_1,n_2n_1n_2,\allowmathbreak n_2n_3n_2,n_3n_2n_3,n_3n_1n_3$ corresponding to these lines, we obtain pairwise
disjoint copies of $T_1$. They intersect $T_1$ only in the listed lines. Denote by $T_2$ the union of $T_1$ with the listed
copies of $T_1$, and~so on. The standard arguments will show that the group $M_3$ acts discretely on $C_1$ and that $T$ is its
fundamental region if we observe that $D:=\bigcup_jT_j$ coincides with $C_1$.

As $D$ is $M_3$-stable and open in $C_1$, its boundary $\partial D:=\overline D\setminus D$ is $M_3$-stable and closed in
$\Bbb R^3$. It~suffices to show that $\partial D=\varnothing$.

Suppose that the intersection $D_r$ of $\partial D$ with the closed ball of radius $r>1$ centred at the origin in $\Bbb R^3$ is
nonempty. Since $D_r$ is compact, there is a point $P\in D_r$ with a minimal value of $t_1t_2t_3$.
As~$T\cap \partial D=\varnothing$, we have, say, $t'_3P<t_3P$. Then $t_1(Pn_3)=t_1P$, $t_2(Pn_3)=t_2P$, $t_3(Pn_3)=t'_3P<t_3P$,
and $Pn_3\in D_r$. This contradicts the choice of $P$.

The remaining case $a_0=a_1=a_2=a_3=\frac\pi2$ is considered in the next remark
$_\blacksquare$

\medskip

{\bf4.6.4.~Remark.} {\sl The connected component\/ $C$ of\/ $S:=S(\frac\pi2,\frac\pi2,\frac\pi2,\frac\pi2)$ with the action of
the group\/ $M$ is isomorphic to the space\/ $\Bbb R^2/\pm1$ with the natural action of the congruence subgroup\/ $\Gamma_2$ in\/
$\SL_2\Bbb Z$ {\rm(}known to be generated by\/ $\left[\smallmatrix1&2\\0&1\endsmallmatrix\right]$ and\/
$\left[\smallmatrix1&0\\2&1\endsmallmatrix\right]${\rm)}. Therefore, almost every\/ $M$-orbit is dense in\/ $C$.}

\medskip

{\bf Proof.} All vertices of an arbitrary hexagon $P\in S$ lie on a geodesic. Therefore, seeing a hexagon as a triple
$(c_1,c_2,c_3)$, we can interpret it as a triple of real numbers considered up to isometries of the real line. In these terms,
$R_j:r\mapsto 2c_j-r$. Since the modifications $n_2n_1$ and $n_3n_1$ transform the relation $R_3R_2R_1R_0=1$ (we place the label
$c_0$ at $p_3$) into the relations $R_3R_2^{R_2R_1}R_1^{R_2}R_0=1$ and $R_3R_2^{R_1R_3}R_1R_0^{R_3R_1}=1$, respectively, assuming
$c_3:=0\in\Bbb R$, we can see that $n_2n_1:(c_1,c_2)\mapsto\big(2c_2-c_1,2c_2-(2c_1-c_2)\big)=(-c_1+2c_2,-2c_1+3c_2)$ and
$n_3n_1:(c_1,c_2)\mapsto(c_1,2c_1+c_2)$. In other words, we identify $S$ with $\Bbb R^2/\pm1$ where $n_2n_1$ acts as
$\left[\smallmatrix-1&2\\-2&3\endsmallmatrix\right]$ and $n_3n_1$, as $\left[\smallmatrix1&0\\2&1\endsmallmatrix\right]$. It
remains to observe that
$\left[\smallmatrix1&0\\1&1\endsmallmatrix\right]\left[\smallmatrix1&2\\0&1\endsmallmatrix\right]\left[\smallmatrix1&0\\-1&1
\endsmallmatrix\right]=\left[\smallmatrix-1&2\\-2&3\endsmallmatrix\right]$
and
$\left[\smallmatrix1&0\\1&1\endsmallmatrix\right]\left[\smallmatrix1&0\\2&1\endsmallmatrix\right]\left[\smallmatrix1&0\\-1&1
\endsmallmatrix\right]=\left[\smallmatrix1&0\\2&1\endsmallmatrix\right]$.

\bigskip

\centerline{\bf5.~Hyperbolic $2$-spheres with $n+1$ cone singularities and convex $2n$-gons}

\medskip

Given $a_0,a_1,\dots,a_n\in[\frac\pi2,\pi[$ such that $\sum_ja_j>2\pi$, we deal with the space $C(a_0,a_1,\dots,a_n)$ of
hyperbolic $2$-spheres $\Sigma$ with labeled cone singularities $c_0,c_1,\dots,c_n$ such that the apex curvature at $c_j$ equals
$2a_j$ for all $j$; the $2$-spheres are considered up to orientation- and label-preserving isometries.

Let $\Sigma\in C$. We pick geodesic segments $[c_0,c_j]\subset\Sigma$ whose pairwise intersections are just $c_0$
(say,~the~shortest ones) and cut $\Sigma$ along these segments. We get a hyperbolic $2n$-gon $P$ with conse\-cutive vertices
$c_1,p_1,\dots,c_n,p_n$ such that the interior angle at $c_j$ equals $2\pi-2a_j\leqslant\pi$ for all $1\leqslant j\leqslant n$
and the interior angles at the $p_j$'s sum to $2\pi-2a_0\leqslant\pi$. It follows that $P$ is star-like, hence, it is embeddable
by Remark 2.11 into the hyperbolic plane as a convex $2n$-gon. Applying to $P$ the modification $m_{j,l}$ (see~Definition 4.3),
we obtain another convex $2n$-gon $Pm_{j,l}$ such that the sphere $\Sigma$ is glued from $Pm_{j,l}$.

In this section, we show that $C(a_0,a_1,\dots,a_n)=S_0(a_0,a_1,\dots,a_n)/M_n$, where $M_n$ is the group generated by the
modifications $m_{j,l}$ and $S(a_0,a_1,\dots,a_n)\supset S_0(a_0,a_1,\dots,a_n)$ is the union of $(n-1)!$ components formed by
all convex $2n$-gons in $S(a_0,a_1,\dots,a_n)$. Components in question correspond to types of $2n$-gons (see Subsection 1.2.5 for
the definition of the type). The other components of $S(a_0,a_1,\dots,a_n)$ are formed by the $2n$-gons that are nonsimple as
closed curves. It is easy to see that such $2n$-gons can also be characterized as those that `limit' area different from
$2\sum_ja_j-(n+4)\pi$ (see~[Ana1] for the definition of the area `limited' by a nonsimple closed curve).

\medskip

{\bf5.1.~Theorem.} {\sl Let\/ $a_0,a_1,\dots,a_n\in[\frac\pi2,\pi[$ be such that\/ $\sum_ja_j>2\pi$. Given convex\/ $2n$-gons\/
$P,P'\in S_0(a_0,a_1,\dots,a_n)$ such that the corresponding\/ $2$-spheres\/ $\Sigma,\Sigma'$ are isometric\/ {\rm(}the
orientation and the labels are preserved\/{\rm)}, there exists\/ $m\in M_n$ such that\/ $P'=Pm$.}

\medskip

{\bf5.2.}~We take a $2n$-gon $P\in S_0(a_0,a_1,\dots,a_n)$ with consecutive vertices $c_1,p_1,\dots,c_n,p_n$. Without loss of
generality, we assume that the interior angles of $P$ at the $p_j$'s sum to $2\pi-2a_0$ and the interior angle of $P$ at $c_j$
equals $2\pi-2a_j$ for all $1\leqslant j\leqslant n$. Denote by $D$ the closed disc in the hyperbolic plane limited by $P$, let
$\Sigma$ stand for the $2$-sphere glued from $P$, and let $p,c_1,\dots,c_n\in\Sigma$ be respectively the cone points of apex
curvatures $2a_0,2a_1,\dots,2a_n$. The sides of $P$ (after gluing) provide simple geodesic segments $s_j\subset\Sigma$ joining
respectively $p$ and $c_j$, $1\leqslant j\leqslant n$, with pairwise intersection $p$.

Any $2n$-gon $P'\in S_0(a_0,a_1,\dots,a_n)$ that produces the same $\Sigma$ after gluing is therefore nothing but simple geodesic
segments $s'_j\subset\Sigma$ joining respectively $p$ and $c_j$, $1\leqslant j\leqslant n$, and $p$ is their pairwise
intersection. Some of the $s_j$'s can coincide with some of the $s'_j$. Denote by $0\leqslant u\leqslant n$ the number of
noncoinciding ones. Clearly, $P'=P$ if $u=0$.

In terms of the closed disc $D$, an arbitrary segment $s'_j$ is given by finitely many disjoint geodesic segments
$[q_0,q'_1],[q_1,q'_2],\dots,[q_{v-1},q'_v]\subset D$, where $1\leqslant v$, $q_0=p_k$ for some $1\leqslant k\leqslant n$,
$q'_v=c_j$, the~interior of $[q_l,q'_{l+1}]$ lives in the interior of $D$ for all $0\leqslant l\leqslant v-1$ unless $v=1$ and
$k=j-1,j$, and, for all $1\leqslant l\leqslant v-1$, the points $q'_l\ne q_l$ are in the interiors of different sides of $P$
(glued in $\Sigma$), symmetric relatively their common vertex $c_{w_l}$, $1\leqslant w_l\leqslant n$, i.e.,
$\dist(q'_l,c_{w_l})=\dist(c_{w_l},q_l)$. Clearly, $s'_j=s_j$ if $v=1$ and $k=j-1,j$.

\medskip

{\bf Proof of Theorem 5.1.} We work in the settings of 5.2 and proceed by induction on $u$ and then, by induction on the number
of intersections of the noncoinciding segments provided by $P'$ with those provided by $P$. As we do not count as intersections
the points $p,c_1,\dots,c_n$, the intersections in question are represented by points in the interior of $D$ or in the interiors
of the sides of $P$.

The case $u=0$ was done in 5.2. Therefore, we assume $u\geqslant1$. This means that there exists $1\leqslant j\leqslant n$ such
that $s'_j\ne s_j$. Hence, we obtain $v\geqslant1$ disjoint geodesic segments as in 5.2.

We claim that the $2n$-gons $Pm_{k,w}$ and $P'$ satisfy the induction hypothesis, where $w:=w_1$ if $v>1$ and $w:=j$ if $v=1$.

Suppose that $v=1$. Since $s'_j\ne s_j$, we conclude that $k\ne j-1,j$ (see 5.2) and that $m_{k,j}$ is indeed applicable.
Moreover,  when we pass from $P$ to $Pm_{k,j}$, the segment $s_j\subset\Sigma$ is substituted by the segment $s'_j\subset\Sigma$
and the rest of the segments provided by $P$ remains the same. In other words, $u$ diminishes.

Suppose that $v>1$. Then $k\ne w_1-1,w_1$ because the interior of $[q_0,q'_1]=[p_k,q'_1]$ lives in the interior of $D$. So,
$m_{k,w}$ is applicable. When passing from $P$ to $Pm_{k,w}$, we simply replace $s_w$ by the segment $s$ represented by
$[c_w,q_0]$. So, $u$ cannot grow and in fact remains the same. Every intersection with $s$ can be seen as an intersection with
$[c_w,q_0]$ living inside the interior of $D$. Looking more closely at the triangle with the vertices $q_0,q'_1,c_w$, we
understand that every intersection of some $s'_l$ with $s$ can be seen as an intersection of $s'\subset D$, one of the disjoint
segments related to $s'_l$ (as in 5.2), with $[c_w,q_0]$. The latter induces an intersection of $s'$ with $[q'_1,c_w]$ because
$s'$ and $[q_0,q'_1]$ are disjoint: they are related to different segments provided by $P'$. It remains to observe that the
intersection corresponding to $q'_1$ disappears when passing from $P$ to $Pm_{k,w}$
$_\blacksquare$

\bigskip

\centerline{\bf References}

\medskip

\noindent
[Ale] A.~D.~Alexandrov, {\sl Convex polyhedra,} Springer Monographs in Mathematics. Springer-Verlag, Berlin, 2005.
Translated from the 1950 Russian edition by N.~S.~Dairbekov, S.~S.~Kutateladze, and A.~B.~Sossinsky, with comments and
bibliography by V.~A.~Zalgaller and appendices by L.~A.~Shor and Yu.~A.~Volkov

\medskip

\noindent
[AGr] S.~Anan$'$in, C.~H.~Grossi, {\sl Coordinate-free classic geometries,} Moscow Math.~J.~{\bf11} (2011), no.~4,
633--655

\medskip

\noindent
[Ana1] S.~Anan$'$in, E.~C.~B.~Goncalves {\sl A hyperelliptic view on Teichm\"uller space. I,} arXiv: 0709.1711

\medskip

\noindent
[Ana2] S.~Anan$'$in, {\sl Reflections, bendings, and pentagons,} arXiv: 1201.1582

\medskip

\noindent
[DM1] P.~Deligne, G.~D.~Mostow, {\sl Monodromy of hypergeometric functions and non-lattice integral monodromy,}
Publ.~Math.~IHES, no.~63 (1986), 5--89

\medskip

\noindent
[DM2] P.~Deligne, G.~D.~Mostow, {\sl Commensurabilities among lattices in $\PU(1,n)$,} Annals of Math.~Studies
{\bf132} (1993)

\medskip

\noindent
[DPP] M.~Deraux, J.~R.~Parker, J.~Paupert, {\sl New non-arithmetic complex hyperbolic lattices,} Invent.~Mat.~{\bf203} (2016),
no.~3, 681--771

\medskip

\noindent
[FaW] E.~Falbel, R.~A.~Wentworth, {\sl On products of isometries of hyperbolic space,} Topology and its Applications {\bf156}
(2009), 2257--2263

\medskip

\noindent
[GhP] S.~Ghazouani, L.~Pirio, {\sl Moduli spaces of flat tori with prescribed holonomy,} arXiv: 1604.01812

\medskip

\noindent
[GMST] W.~M.~Goldman, G.~McShane, G.~Stantchev, S.~P.~Tan, {\sl Automorphisms of two-generator free groups and spaces of
isometric actions on the hyperbolic plane,} arXiv: 1509.03790v2

\medskip

\noindent
[Gol] W.~M.~Goldman, {\sl The modular group action on real\/ $\SL(2)$-characters of a one-holed torus,}
Geometry \& Topology {\bf7} (2003), 443--486

\medskip

\noindent
[HRi] C.~D.~Hodgson, I.~Rivin, {\sl A characterization of compact convex polyhedra in hyperbolic\/ $3$-space,}
Invent.~Mat.~{\bf111} (1993), 77--111

\medskip

\noindent
[Kui] N.~H.~Kuiper, {\sl Hyperbolic\/ $4$-manifolds and tessellations,} Publ.~Math.~IHES, no.~68 (1988), 47--76

\medskip

\noindent
[McR] D.~B.~McReynolds, {\sl Arithmetic lattices in\/ $\SU(n,1)$,} 2015,\newline
http://www.its.caltech.edu/$\widetilde{\phantom{m}}$dmcreyn/ComplexArithmeticI.pdf

\medskip

\noindent
[MPa] G.~Mondello, D.~Panov, {\sl Spherical metrics with conical singularities on a\/ $2$-sphere\/{\rm:} angle constraints,}
arXiv: 1505.01994

\medskip

\noindent
[Thu] W.~P.~Thurston, {\sl Shapes of polyhedra and triangulations of the sphere,} Geometry \& Topology Monographs
{\bf1} (1998), 511--549

\medskip

\noindent
[WaW] P.~Waterman, S. Wolpert, {\sl Earthquakes and tessellations of Teichm\"uller space,} TAMS {\bf278} (1983), no.~1,
157--167

\enddocument